\documentclass[12pt]{article}
\usepackage{amscd,amsfonts,amssymb,amsmath,latexsym,array,hhline}
\usepackage[dvips]{graphics}
\makeatletter
\@addtoreset{equation}{section}
\renewcommand{\@begintheorem}[2]{\begin{trivlist}\it
\item[\hspace{\labelsep}{\bf #1\ #2.}]}
\renewcommand{\@opargbegintheorem}[3]{\begin{trivlist}\it
\item[\hspace{\labelsep}{\bf #1\ #2\ (#3).}]}
\renewcommand{\@endtheorem}{\end{trivlist}}
\renewcommand{\@cite}[2]{[{#1\if@tempswa ; #2\fi}]}
\makeatother

\newcommand{\paragr}{\hspace{6mm}}
\newtheorem{Theorem}{\paragr Theorem}[section]
\newtheorem{Lemma}[Theorem]{\paragr Lemma}

\newtheorem{Proposition}[Theorem]{\paragr Proposition}
\newtheorem{Definition}[Theorem]{\paragr Definition}

\newtheorem{Example}[Theorem]{\paragr Example}
\newtheorem{Corollary}[Theorem]{\paragr Corollary}
\newcommand{\Proof}{\texttt{Proof}. }
\newcommand{\Remark}{\texttt{Remark}. }
\newcommand{\Iremark}{\texttt{Important remark}. }

\newcommand{\Remarks}{\texttt{Remarks}. }

\newcommand{\End}{~\hfill $\Box\!$}

\newcommand{\sks}{\smallskip}
\newcommand{\skm}{\medskip}
\newcommand{\skb}{\bigskip}
\newcommand{\sksminus}{}
\newcommand{\sksmin}{}

\newcommand{\al}{\alpha}

\newcommand{\ga}{\gamma}

\newcommand{\la}{\lambda}
\newcommand{\si}{\sigma}

\newcommand{\eps}{\varepsilon}
\renewcommand{\phi}{\varphi}
\renewcommand{\kappa}{\varkappa}

\newcommand{\R}{\mathbb{R}}

\newcommand{\F}{\mathcal{F}}

\newcommand{\DDD}{\mathcal{D}}
\newcommand{\EEE}{\mathcal{E}}

\newcommand{\RRR}{\mathcal{R}}
\newcommand{\PPP}{\mathcal{P}}
\newcommand{\EX}{\mathcal{X}}
\newcommand{\PD}{\mathcal{PD}}
\newcommand{\RD}{\mathcal{RD}}
\newcommand{\EE}{\mathsf{E}}
\newcommand{\PP}{\mathsf{P}}
\newcommand{\QQ}{\mathsf{Q}}

\newcommand{\emp}{\emptyset}
\newcommand{\lb}{\langle}
\newcommand{\rb}{\rangle}
\newcommand{\wt}{\widetilde}
\newcommand{\wl}{\overline}
\newcommand{\xra}{\xrightarrow}
\newcommand{\da}{\downarrow}

\newcommand{\ds}{\displaystyle}

\newcommand{\conv}{\mathop{\rm conv}\nolimits}
\newcommand{\cone}{\mathop{\rm cone}}

\newcommand{\pr}{\mathop{\rm pr}\nolimits}

\newcommand{\argmax}{\mathop{\rm argmax}}
\newcommand{\argmin}{\mathop{\rm argmin}}

\newcommand{\cl}{\mathop{\rm cl}}
\newcommand{\cov}{\mathop{\sf cov}}
\newcommand{\var}{\mathop{\sf var}}

\renewcommand{\inf}{\mathop{\rm inf\rule[-0.8mm]{0mm}{1mm}}}

\newcommand{\RAROC}{{\rm RAROC}}

\newcommand{\INGDR}{I_{\text{\sl NGD}\,(\!R)}}

\newcommand{\INBC}{I_{\text{\sl NBC}}}
\newcommand{\IE}{I_{\text{\sl E}}}
\newcommand{\St}{H}

\newenvironment{mitemize}%
{\begin{list}{$\bullet$}{
\leftmargin=32pt
\rightmargin=0pt
\labelsep=5pt
\labelwidth=20pt
\itemindent=0pt
\topsep=5pt plus 2pt minus 4pt
\partopsep=2pt plus 1pt minus 1pt
\parsep=0pt
\itemsep=0pt}}%
{\end{list}}

\begin{document}
\vspace*{1mm}
\begin{center}\bf
EQUILIBRIUM WITH COHERENT RISK
\end{center}

\begin{center}\itshape\bfseries
A.S.~Cherny
\end{center}

\begin{center}
\textit{Moscow State University,}\\
\textit{Faculty of Mechanics and Mathematics,}\\
\textit{Department of Probability Theory,}\\
\textit{119992 Moscow, Russia.}\\
\texttt{E-mail: cherny@mech.math.msu.su}\\
\texttt{Webpage: http://mech.math.msu.su/\~{}cherny}
\end{center}

\begin{abstract}
\textbf{Abstract.}
This paper is the continuation of~\cite{C061}
and deals with further applications of coherent risk
measures to problems of finance.

First, we study the optimization problem.
Three forms of this problem are considered.

Furthermore, the results obtained are
applied to the optimality pricing.
Again three forms of this technique are considered.

Finally, we study the equilibrium problem
both in the unconstrained and in the constrained forms.
We establish the equivalence between the global and the
competitive optima and give a dual description of the
equilibrium. Moreover, we provide an explicit geometric
solution of the constrained equilibrium problem.

Most of the results are presented on two levels:
on a general level the results have a probabilistic form;
for a static model with a finite number of assets,
the results have a geometric form.

\end{abstract}

%======================================================
\section{Introduction}
\label{I}

\textbf{1. Goal of the paper.}
In this paper, we consider applications of coherent
risk measures~to
\begin{mitemize}
\item optimization;
\item optimality pricing;
\item equilibrium.
\end{mitemize}

The optimization problem is considered in three forms.
First we study what we call the
\textit{agent-independent optimization}.
It is in fact the Markowitz-type optimization problem
with variance replaced by a coherent risk
measure,\footnote{It has been clear from the outset
that variance is not a very good measure of risk
because high profits are penalized in the same way as
high losses. In~\cite{M59}, Markowitz proposed a way
to overcome this problem by considering semivariance
$\|(X-\EE X)^-\|_{L^2}$ instead of variance.
The function $\rho(X)=-\EE X+\al\|(X-\EE X)^-\|_{L^2}$
with $0\le\al\le1$ is, in fact, an example of a coherent risk
measure (see~\cite{F03}). Thus, in essence, semivariance
is a particular case of the coherent risk.}
i.e. a problem of the form
\begin{equation}
\label{i1}
\begin{cases}
\EE X\longrightarrow\max,\\
\rho(X)\le c,
\end{cases}
\end{equation}
where $X$ means the discounted P\&L earned by a portfolio
and $\rho$ is a coherent risk measure
(P\&L means the Profit\&Loss, i.e. the difference
between the terminal wealth and the initial wealth).
Let us remark that this problem was considered
in~\cite{A04}, \cite{RU00}, \cite{RUZ05}.\footnote{If~$\rho$
is defined by a finite number of probabilistic scenarios,
then the corresponding optimization problem becomes the one
considered in the generalized Neyman--Pearson lemma
(see~\cite[Ch.~3]{L97}).}
As opposed to these papers, we have at our disposal the
notion of a \textit{generator} introduced in~\cite{C061}.
In terms of generators, we are able to give a geometric
solution of~\eqref{i1} (see Figures~1, 2).
The model we are considering takes into account such
market imperfections as cone portfolio constraints,
transaction costs, and the ambiguity of the historic
probability measure.

Problem~\eqref{i1} is the optimization problem for an
investor whose capital evolves in a risk-free way.
However, an investor might have a risky endowment
with a random terminal wealth~$W$
($W$ might have a financial or a non-financial
structure; for example, it might be the terminal
wealth of a firm producing some goods).
The investor trying to minimize his/her risk by trading
in the market faces the problem
\begin{equation}
\label{i2}
\rho(X+W)\longrightarrow\min.
\end{equation}
Note that this coincides with the superreplication problem
for the NGD pricing (see~\cite[Subsect.~3.6]{C061}).
We study~\eqref{i2} in two forms, which we call the
\textit{global} and the \textit{local single-agent optimization},
respectively. On the financial side, the former pertains
to a ``small'' investor, while the latter pertains to a
``big'' investor.
For both of them, we provide a geometric solution
(see Figures~4, 6).
Also, in~\cite[Sect.~5]{C05e}, we provide sufficient
conditions for the uniqueness of a solution of various
optimization problems.

After considering the optimization problem, we return
to the pricing problem and apply the obtained
results to the optimality pricing of contingent
claims (in the finance literature, this is typically
referred to as reservation pricing).
Again, we propose three forms of this technique.
They differ by the inputs they require and the assumptions
behind them (thus, of course, they produce different
outputs).

The first technique, which we call the \textit{agent-independent
optimality pricing}, might be considered as the limit
case of the RAROC-based NGD pricing technique introduced
in~\cite{C061}.
A rather surprising outcome of this technique is that
it provides (typically) a single fair price of a
contingent claim for all agents (under the strong assumption
that all the agents are using the same historic measure~$\PP$
and the same risk measure~$\rho$ and
all the agents are trying to maximize
RAROC defined through~$\PP$ and~$\rho$).

Another technique is called the \textit{single-agent
optimality pricing}. It is, in fact, the coherent version
of the classical utility-indifference pricing with the
expected utility replaced by the coherent one.
As an outcome, this technique typically provides a single
number, which means the price of a contingent claim that
is fair for a particular agent (it depends on the risk
measure he/she is using and on his/her endowment).

Then we consider one more technique called the
\textit{multi-agent optimality pricing}.
The idea is as follows. We have a contingent claim and
several agents, each employing his/her own coherent risk
measure. A price is said to be fair if it provides no
trading opportunity, which would allow each agent to
decrease his/her risk (this is a modified form
of the idea proposed in~\cite{CGM01}).
As an outcome, this technique typically provides a whole
interval of prices, which are fair for this group of
agents (in fact, this interval is the convex hull of
fair prices for these agents produced by the previous
technique).

Finally, we turn to the equilibrium problem.
One of the basic results of the classical economic theory
is the equivalence between the Pareto optimum
(known also as the ``Soviet-type optimum'') and the
competitive optimum (known also as the ``western-type
optimum''); see~\cite{K90}.
This result is established within the framework of the
expected utility.

In the present paper, we establish the analog of this result
within the coherent utility framework.
This is done for two types of equilibrium:
for the unconstrained one (Theorem~\ref{EC7})
and for the constrained one (Theorem~\ref{EI5}).
A very important feature of coherent utility (which is
not shared by the expected utility) is that it admits a
rich duality theory. Thus, we not only establish the equivalence
between different types of equilibrium, but also provide
its dual description.

Moreover, for the constrained equilibrium problem, we
are able to provide an explicit geometric solution
based on generators (see Figure~11).
It yields the equilibrium price as well as the equilibrium
portfolios of the agents.

\skb
\textbf{2. Structure of the paper.}
Section~\ref{O} deals with the optimization problem.
In Subsections~\ref{OG}, \ref{GOG}, and~\ref{LOG},
we consider the three techniques described above.
In Subsection~\ref{OS}, the obtained results are applied
to the problem of finding the optimal structure of a firm
consisting of several units.
The provided theorem states that a structure is optimal
if and only if the \textit{RAROC contributions} of
different units are the same.
In Subsection~\ref{UPL}, we apply the obtained results
to the study of the liquidity effects in the framework
of the NGD pricing considered in~\cite{C061}.

Section~\ref{OP} is related to the optimality pricing.
Subsections~\ref{AIOP}, \ref{SAOP}, and~\ref{MAOP}
correspond to the three techniques described above.

Section~\ref{E} deals with equilibrium.
Subsections~\ref{EC} and~\ref{EI} are, in fact, duals
of each other: they consider the unconstrained and the
constrained equilibria, respectively.

Altogether, there are seven pricing techniques proposed
in~\cite{C061} and in the present paper.
They are compared in the final Section~\ref{C}.

%=======================================================
\section{Optimization}
\label{O}

%-------------------------------------------------------
\subsection{Agent-Independent Optimization}
\label{OG}

We consider the model of~\cite[Subsect.~3.2]{C061}.
Thus, we are given a probability space $(\Omega,\F,\PP)$,
a convex weakly compact set $\RD\subset\PPP$,
an $L^1$-closed convex set $\PD\subset\RD$,
and a convex set $A\subset L^0$.
Let us introduce the notation
$\EE_{\PD}X=\inf_{\QQ\in\PD}\EE_\QQ X$,
$u(X)=\inf_{\QQ\in\RD}\EE_\QQ X$,
$\rho(X)=-u(X)$ (we understand $\EE_\QQ X$ according to
the convention of~\cite[Def.~2.3]{C061}).

\skm
\textbf{Problem (agent-independent optimization):}
The problem is
$$
\begin{cases}
\EE_{\PD}X\longrightarrow\max,\\
X\in A,\;\rho(X)\le c,
\end{cases}
$$
where $c\in\R_+$.
Clearly, if $A$ is a cone, then this problem is obviously
equivalent to the problem of finding
$$
R_*=\sup_{X\in A}\RAROC(X)
$$
and
$$
X_*=\argmax_{X\in A}\RAROC(X),
$$
where
$$
\RAROC(X)=\begin{cases}
+\infty&\text{if}\;\;\EE_\PD X>0
\text{ and }u(X)\ge0,\\[2mm]
\ds\frac{\EE_\PD X}{\rho(X)}
&\text{otherwise}
\end{cases}
$$
with the convention $\frac{0}{0}=0$,
$\frac{\infty}{\infty}=0$.

The only statement we can make at this level of generality
is that
$$
R_*=\inf\biggl\{R>0:\Bigl(\frac{1}{1+R}\,\PD
+\frac{R}{1+R}\,\RD\Bigr)\cap\RRR\ne\emp\biggr\},
$$
which follows from~\cite[Th.~3.10]{C061}.
Of course, in general $X_*$ need not exist.

We will now study the problem for a static model with
a finite number of assets.
Let $A=\{\lb h,S_1-S_0\rb:h\in\St\}$, where
$S_0\in\R^d$, $S_1^1,\dots,S_1^d\in L_w^1(\RD)$, and
$\St\subseteq\R^d$ is a closed convex cone
(here we impose no conditions on~$\DDD$).
Let us introduce the notation (see Figure~1)
\begin{equation}
\begin{split}
\label{of1}
\St^*&=\{x\in\R^d:\forall h\in\St,\:\lb h,x\rb\ge0\},\\
E&=\cl\{\EE_\QQ S_1:\QQ\in\PD\},\\
G&=\cl\{\EE_\QQ S_1:\QQ\in\RD\},\\
D&=G+\St^*,
\end{split}
\end{equation}
where ``$\cl$'' denotes the closure,
and let $D^\circ$ denote the relative interior of~$D$.
(The set $G$ is the generator for $S_1$ and $u$.)
The sets $E$ and $G$ are convex compacts,
while $D$ is convex and closed.
Note that, for $h\in\St$,
\begin{align}
\label{of2}
\EE_\PD\lb h,S_1-S_0\rb
&=\inf_{x\in E}\lb h,x-S_0\rb,\\
\label{of3}
u(\lb h,S_1-S_0\rb)
&=\inf_{x\in G}\lb h,x-S_0\rb
=\inf_{x\in D}\lb h,x-S_0\rb.
\end{align}

We will assume that $S_0\in D^\circ\setminus E$.
This assumption is justified economically. Indeed,
if $S_0\in E$, then, in view of~\eqref{of2},
$\RAROC(X)=0$ for any $X\in A$;
if $S_0\notin D^\circ$, then, in view of~\eqref{of3},
there exists $X\in A$ with $\RAROC(X)=\infty$
(provided that $E$ belongs to the relative interior of~$C$).

For $\la>0$, we denote $E(\la)=S_0-\la(E-S_0)$
and set $\la_*=\sup\{\la>0:E(\la)\cap D\ne\emp\}$,
$$
N=\{h\!\in\!\St:\exists a\!\in\!\R:
\forall x\!\in\! E(\la_*),\forall y\!\in\! D,\;
\lb h,x\rb\le a\le\lb h,y\rb
\text{ and }
\forall y\!\in\! D^\circ\!,\;
\lb h,y\rb>a\}.
$$
Note that $N$ is nonempty provided that $\la_*<\infty$.
In the case, where $\la_*=\infty$, we set $N=\St$.

\begin{figure}[!h]
\begin{picture}(150,62)(-31,-19.5)
\put(0,0){\includegraphics{poem.5}}
\put(74,12.5){\small $E$}
\put(71,24){\small $G$}
\put(70.5,34){\small $D$}
\put(43.5,20.7){\scalebox{0.75}{$E(\!\la_*\!)$}}
\put(28,24.5){\small $E(\la)$}
\put(58,13.5){\small $S_0$}
\put(62.5,23.5){\small $h_*$}
\put(17.5,6){\small $\St$}
\put(10,16){\small $\St^*$}
\put(15,-15){\parbox{67mm}{\small\textbf{Figure~1.}
Solution of the optimization problem.
Here $h_*$ is an optimal $h$.}}
\end{picture}
\end{figure}

\begin{Theorem}
\label{OF1}
We have $R_*=\la_*^{-1}$ and
$\argmax_{h\in\St}\RAROC(\lb h,S_1-S_0\rb)=N$.
\end{Theorem}

\sksminus
\Proof
We will prove the statement for the case $\la_*<\infty$.
The proof for the case $\la_*=\infty$ is similar.
Take $T\in E(\la_*)\cap D$ and set $U=S_0-\la_*^{-1}(T-S_0)$.

If $h\in N$, then
$$
\RAROC(\lb h,S_1-S_0\rb)
=\frac{\inf_{x\in E}\lb h,x-S_0\rb}{-\inf_{x\in D}\lb h,x-S_0\rb}
=\frac{\lb h,U-S_0\rb}{-\lb h,T-S_0\rb}
=\la_*^{-1}.
$$

If $h\in\St\setminus N$, then there are three possibilities:
\begin{mitemize}
\item[1)] $h$ is orthogonal to the smallest affine subspace
containing~$D$;
\item[2)] $\sup_{x\in E(\la_*)}\lb h,x\rb>\lb h,T\rb$;
\item[3)] $\inf_{x\in D}\lb h,x\rb<\lb h,T\rb$.
\end{mitemize}
In the first case, $\RAROC(\lb h,S_1-S_0\rb)=0$.
In the second case,
$$
\inf_{x\in E}\lb h,x-S_0\rb<\lb h,U-S_0\rb,\qquad
\inf_{x\in D}\lb h,x-S_0\rb\le\lb h,T-S_0\rb,
$$
so that $\RAROC(\lb h,S_1-S_0\rb)<\la_*^{-1}$.
The third case is analyzed in a similar way.\End

\skm
As a corollary, in the case, where $\PD=\{\PP\}$ and
$\St=\R^d$, the solution to the optimization problem
is found as follows.
Let~$T$ be the intersection of the ray
$(E,S_0)$ (in this case $E=\EE_\PP S_1$) with the
border of~$G$. Then
$$
\sup_{h\in\R^d}\RAROC(\lb h,S_1-S_0\rb)
=\frac{|\EE_\PP S_1-S_0|}{|S_0-T|}
$$
and $\argmax_{h\in\R^d}\RAROC(\lb h,S_1-S_0\rb)$ is
$$
N_G(T):=\{h\in\R^d:\forall x\in G^\circ,\;\lb h,x-T\rb>0\}.
$$
In the case, where $G$ has a nonempty interior, $N_G(T)$
is the set of inner normals to $G$ at the point~$T$.

\begin{figure}[!h]
\begin{picture}(150,52.5)(-80,-35)
\put(-28.5,-16){\includegraphics{poem.6}}
\put(-22,-14){\small $T$}
\put(-3,-5){\small $S_0$}
\put(4,-1.5){\small $E$}
\put(-14,0){\small $h_*$}
\put(14,-1){\small $G$}
\put(-38,-30){\parbox{75mm}{\small\textbf{Figure~2.}
Solution of the optimization problem in the case
$\PD=\{\PP\}$ and $H=\R^d$}}
\end{picture}
\end{figure}

\skm
\Iremark
In order to find the solution of the optimization
problem for the case $\PD=\{\PP\}$, $\St=\R^d$, one
needs to know the generating set $G$ and the vector
$\EE_\PP S_1$. The empirical estimation of $G$ is a
problem similar to the empirical estimation of
volatility, and hence, it can be successfully accomplished.
However, the empirical estimation of the mean vector
$\EE_\PP S_1$ is known to be a very unpleasant problem
because it is very close to~0
(see the discussion in~\cite{B95} and
the 20's example in~\cite{JPR05}).
But it turns out that the well-known security market line
relationship of Sharpe~\cite{S64} helps to overcome this
problem. This relation states that
$$
\EE_\PP\biggl(\frac{\wt S_1^i-\wt S_0^i}{\wt S_0^i}-r\biggr)
=\beta^i\EE_\PP\biggl(\frac{\wt S_1^M-\wt S_0^M}{\wt
S_0^M}-r\biggr),\quad i=1,\dots,d,
$$
where $r$ is the risk-free interest rate,
$\wt S_n^i=(1+r)^nS_n^i$ are true (not discounted) prices,
and $\wt S_n^M$ is the price of the market portfolio at time~$n$.
Hence,
$$
\EE_\PP(S_1^i-S_0^i)=\beta^i const,\quad i=1,\dots,d.
$$
The constant here contains as a factor the expected
excess return on the
market portfolio, which is again hard to estimate.
But note that for our purposes this unknown constant
is not needed!
Indeed, the geometric solution of the optimization problem
presented above requires only the direction of the vector
$\EE_\PP S_1-S_0$, and this depends only on
$(\beta^1,\dots,\beta^d)$.

\skm
The following example shows that in natural situations
the set of optimal strategies~$h_*$ might not be unique
(of course, the uniqueness of~$h_*$ should be understood
up to multiplication by a positive constant).

\begin{Example}\rm
\label{OF2}
Let $S_1^1$ have a continuous distribution with
$\EE S_1<\infty$ and take $S_1^2=(S_1^1-K)^+$
(so that the second asset is a call option on the first one).
Let $\PD=\{\PP\}$, $\RD$ be the determining set
of $\text{Tail V@R}$ of order $\la$
(see~\cite[Ex.~2.5]{C061}) and $\St=\R^2$.
Assume that $\F=\si(S_1^1)$.
It is easy to see that $\EX_\RD(S_1^1)$ consists of a
unique element $\QQ=\la^{-1}I(S_1^1\le q_\la)\PP$,
where $q_\la$ is the $\la$-quantile of~$S_1^1$.
The border of~$G$ has an angle $\pi/4$
at the point $T=\EE_\QQ(S_1^1,S_1^2)$ (see Figure~3).
Let $S_0=\frac{T+\EE_\PP S_1}{2}$.
Then $N_G(T)=\{h\in\R^2:h^1\ge0,\;h^2\ge-h^1\}$.\End
\end{Example}

\begin{figure}[h]
\begin{picture}(150,55)(-45,-17.5)
\put(-0.5,-7.3){\includegraphics{poem.7}}
\put(0,0){\vector(1,0){70}}
\put(0,0){\vector(0,1){35}}
\put(45.3,16.5){\small $G$}
\put(0,1.5){\scalebox{0.7}{$N_G(T)$}}
\put(13,-4){\small $T$}
\put(26.5,3){\small $S_0$}
\put(39,10){\small $E$}
\put(33,-4){\small $K$}
\put(-6.5,-15){\small\textbf{Figure~3.}
Nonuniqueness of an optimal strategy}
\end{picture}
\end{figure}

Let us now find the solution of the optimization problem
in the Gaussian case.

\begin{Example}\rm
\label{OF3}
Let $S_1$ have Gaussian distribution with mean $a$ and
covariance matrix~$C$.
Let $\PD=\{\PP\}$, $\St=\R^d$, and $\RD$ be the
determining set of a law invariant coherent utility
function~$u$ that is finite on Gaussian random variables.
Assume that $S_0$ belongs to the relative interior of~$G$
and $S_0\ne a$.

There exists $\ga>0$ such that, for a Gaussian random
variable $\xi$ with mean~$m$ and variance~$\si^2$, we have
$u(\xi)=m-\ga\si$.
Let $L$ denote the image of $\R^d$ under the map
$x\mapsto Cx$.
It is easy to see that
$$
G=a+\{C^{1/2}x:\|x\|\le\ga\}
=a+\{y\in L:\lb y,C^{-1}y\rb\le\ga^2\}.
$$
We have $T=a+\al(S_0-a)$ with some $\al>0$.
It is easy to see that $h\in N_D(T)$ if and only if
$\lb h,a-S_0\rb>0$ and, for any $y\in L$ such that
$$
\frac{d}{d\eps}\Bigl|_{\eps=0}
\lb T-a+\eps y,C^{-1}(T-a+\eps y)\rb=0,
$$
we have $\lb\pr_L\!h,y\rb=0$.
This means that
$\pr_L\!h=c' C^{-1}(a-T)=c C^{-1}(a-S_0)$ with some
constant $c>0$.
Thus,
$$
N_D(T)=\{h\in\R^d:Ch=c(a-S_0),\;c>0\}.
$$
Note that this set does not depend on~$u$!

It is easy to see that
$$
R_*
=\frac{|S_0-a|}{|T-S_0|}
=\frac{|S_0-a|}{|T-a|-|S_0-a|}
=\frac{\lb S_0-a,C^{-1}(S_0-a)\rb^{1/2}}%
{\ga-\lb S_0-a,C^{-1}(S_0-a)\rb^{1/2}}.
$$
This equality can also be deduced from~\cite[Ex.~3.13]{C061}.\End
\end{Example}

%--------------------------------------------------------
\subsection{Optimal Structure of a Firm}
\label{OS}

Let $(\Omega,\F,\PP)$ be a probability space,
$\DDD\subseteq\PPP$ be a convex set
(we assume that $\PP\in\DDD$)
and let $X^1,\dots,X^d\in L_w^1(\DDD)$ be the discounted
P\&Ls produced by different components of some firm.

We will consider the problem
\begin{equation}
\label{os1}
\begin{cases}
\EE_\PP\lb h,X\rb\longrightarrow\max,\\
h\in\R_+^d,\;\rho(\lb h,X\rb)\le c,
\end{cases}
\end{equation}
where $c$ is a positive constant meaning the capital
available to the whole firm.
From the financial point of view, \eqref{os1}
is the problem of the central management of the firm
deciding which components should grow and which should
shrink.

This is a particular case of the optimization
problem of the previous subsection
(with $\PD=\{\PP\}$, $\RD=\DDD$, and $\St=\R_+^d$),
so that we already have a geometric recipe to find the
optimal solution. Here we will present an economic
characterization of optimality.
We will consider an arbitrary convex cone constraint~$\St$
(not only $\R_+^d$ as in~\eqref{os1}).
We assume that $\EE_\PP X\ne0$ and that the
generator~$G$ given by~\eqref{of1} is strictly convex,
i.e. its interior is nonempty and its border contains no interval.

\begin{Definition}\rm
\label{OS1}
We define the \textit{RAROC contribution of $X$ to $Y$} as
$$
\RAROC^c(X;Y)
=\frac{\EE_\PP X}{\rho^c(X;Y)},
$$
where $\rho^c$ is the risk contribution
(see~\cite[Subsect.~2.5]{C061}).
\end{Definition}

The RAROC contribution is well defined provided that
$\rho^c(X;Y)$ is well defined and $\rho^c(X;Y)\ne0$.

\skm
\Remarks
\texttt{(i)} The RAROC contribution may take on negative values.

\texttt{(ii)} We have $\RAROC^c(X;X)=\RAROC(X)$.

\begin{Theorem}
\label{OS2}
If $h\in H$ and
\begin{equation}
\label{os2}
\RAROC^c\Bigl(h^1X^1;\sum_{i=1}^d h^iX^i\Bigr)
=\dots=
\RAROC^c\Bigl(h^dX^d;\sum_{i=1}^d h^iX^i\Bigr),
\end{equation}
then $h\in\argmax_{h\in\St}\RAROC(\lb h,X\rb)$
and all the elements of this equality
are equal to~$R_*$.

Conversely, if $h$ is an inner point of~$\St$
and $h\in\argmax_{h\in\St}\RAROC(\lb h,X\rb)$,
then~\eqref{os2} is satisfied.
\end{Theorem}

\Proof
Denote $\sum h^iX^i$ by~$Y$.
It is seen from~\cite[Th.~2.16]{C061} that
$u^c(h^iX^i;Y)=h^i u^c(X^i;Y)$.
Repeating the arguments of the proof of~\cite[Th.~2.12]{C061},
we get $u^c(X^i;Y)=U^i$, where $U=\argmin_{x\in G}\lb h,x\rb$
(this point is unique due to the convexity of~$G$).
Thus, \eqref{os2} is equivalent to: $\EE_\PP X=-R U$,
where $R=\RAROC^c(h^iX^i;Y)$.
It is seen from the results of the previous subsection
that this condition implies that
$h\in\argmax_{h\in\R^d}\RAROC(\lb h,X\rb)$.
As $u(\lb h,X\rb)=\lb h,U\rb$, we get
$\RAROC(\lb h,X\rb)=R$, so that $R_*=R$.

Conversely, if $h$ is an inner point of~$\St$, then
$\argmin_{x\in D}\lb h,x\rb=\argmin_{x\in G}\lb h,x\rb=U$
($D$ is given by~\eqref{of1}).
Recalling the results of the previous subsection, we get
the second statement.\End

\skm
\Remark
The additional assumption that $h$ is in the interior
of~$\St$ is essential for the converse statement of
Theorem~\ref{OS2}.
As an example, take $\St=\{\al h_0:\al\in\R_+\}$,
where $h_0$ is a fixed vector. Then clearly
$h_0\in\argmax_{h\in\St}\RAROC(\lb h,X\rb)$, but of
course~\eqref{os2} might be violated.

%--------------------------------------------------------
\subsection{Single-Agent Global Optimization}
\label{GOG}

Let $(\Omega,\F,\PP)$ be a probability space,
$u$ be a coherent utility function with the weakly compact
determining set~$\DDD$,
$A\subseteq L^0$ be a $\DDD$-consistent convex cone,
and $W\in L_s^1(\DDD)$.
From the financial point of view, $W$ is the
terminal endowment of some agent,
while $A$ is the set of discounted P\&Ls
the agent can obtain by trading.

\skm
\textbf{Problem (single-agent global optimization):}
Find
$$
u_*=\sup_{X\in A}u(W+X)
$$
and
$$
X_*=\argmax_{X\in A}u(W+X).
$$

\begin{Proposition}
\label{GOG1}
We have
$$
u_*=\inf_{\QQ\in\DDD\cap\RRR}\EE_\QQ W,
$$
where $\inf\emp:=\infty$.
\end{Proposition}

\sksmin
\Proof
By~\cite[Th.~3.4]{C061}, for any $z\in\R$,
$$
\sup_{X\in A}u(W+X)>z
\;\Longleftrightarrow\;
\sup_{X\in A}u(-z+W+X)>0
\;\Longleftrightarrow\;
\DDD\cap\RRR(-z+W+A)=\emp
$$
(the notation $\RRR(A)$ was introduced in~\cite[Def.~3.1]{C061}).
Fix $\QQ\in\RRR(-z+W+A)$. As $A$ is a cone, we have
$\EE_\QQ X\le0$ for any $X\in A$.
As $A$ contains zero, $\EE_\QQ(-z+W)\le0$.
Thus, $\QQ\in\RRR$ and $\EE_\QQ(-z+W)\le0$.
Conversely, if $\QQ\in\DDD$ and these two conditions
are satisfied, then $\QQ\in\RRR(-z+W+A)$.
We get
$$
\sup_{X\in A}u(W+X)>z
\;\Longleftrightarrow\;
\DDD\cap\RRR\cap\{\QQ:\EE_\QQ W\le z\}=\emp,
$$
and the result follows.\End

\skm
We will now study the problem for a static model with
a finite number of assets.
Let $A=\{\lb h,X\rb:h\in\St\}$, where
$X=(X^1,\dots,X^d)\in L_w^1(\DDD)$
and $\St\subseteq\R^d$ is a closed convex cone
(here we impose no conditions on~$\DDD$).
For the case $H=\R^d$, we provided a geometric solution
of this problem in~\cite[Subsect.~3.6]{C061}. For an
arbitrary~$H$, it is more complicated and is given below.
Let us introduce the notation (see Figure~4)
\begin{align*}
G&=\cl\{\EE_\QQ(X,W):\QQ\in\DDD\},\\
\wt\St&=\{x\in\R^{d+1}:(x^1,\dots,x^d)\in\St,\;x^{d+1}=1\},\\
\wt\St^*&=\{x\in\R^{d+1}:\forall h\in\wt\St,\;\lb h,x\rb\le0\},\\
e&=(0,\dots,0,1),\\
\la_*&=\inf\{\la\in\R:(\la e+\wt\St^*)\cap G\ne\emp\},\\
\wt N&=\{h\in\R^{d+1}:h^{d+1}=1\text{ and }\exists a\in\R:\\
&\hspace*{7mm}\forall x\in\la_*e+\wt\St^*,\:\forall y\in G,\;
\lb h,x\rb\le a\le\lb h,y\rb\},\\
N&=\{h\in\R^d:(h,1)\in\wt N\}.
\end{align*}
If $\la_*=\infty$, we set $\wt N=N=\emp$.

\begin{figure}[h]
\begin{picture}(150,90)(-67,-37.5)
\put(-20.4,-16.2){\includegraphics{poem.8}}
\put(0,0){\vector(1,0){45}}
\put(0,0){\vector(0,1){50}}
\multiput(9,0.5)(0,2){20}{\line(0,1){1}}
\multiput(0,40)(2,0){5}{\line(1,0){1}}
\put(-3,34){\scalebox{1}{$\biggl\{$}}
\put(-4.5,34){\scalebox{0.9}{$1$}}
\put(43,-4){\small $\R^d$}
\put(-4,48){\small $\R$}
\put(-3,8){\small $e$}
\put(7,-4){\small $h_*$}
\put(9,38.5){\small $\wt h_*$}
\put(15.5,34.8){\small $G$}
\put(11,11){\small $\wt\St$}
\put(-3,-13.6){\small $\wt\St^*$}
\put(0,16){\scalebox{0.85}{\small $\wt\St_\la^*$}}
\put(-21,-30){\parbox{68mm}{\small\textbf{Figure~4.}
Solution of the optimization problem.
By $\wt\St_\la^*$ we denote $\la_*e+\wt\St^*$.
Here $\wt N=\{\wt h_*\}$ and $N=\{h_*\}$.}}
\end{picture}
\end{figure}

\begin{Theorem}
\label{GOF1}
We have $u_*=\la_*$ and
$\argmax_{h\in\St}u(W+\lb h,X\rb)=N$.
\end{Theorem}

\sksminus
\Proof
Fix $\la<\la_*$.
As $G$ is a convex compact and $\wt\St^*$ is convex and closed, there exist
$\wt h\in\R^{d+1}$ and $a,b\in\R$ such that,
for any $x\in\la e+\wt\St^*$ and any $y\in G$, we have
$\lb\wt h,x\rb\le a<b\le\lb\wt h,y\rb$.
As $G$ is compact, $\wt h$ can be chosen in such a way that
$\wt h^{d+1}\ne0$.
Since $\wt\St^*\supseteq\{\al e:\al\le0\}$, we have
$\wt h^{d+1}>0$.
Without loss of generality, $\wt h^{d+1}=1$.
Then, for any $x\in\wt\St^*$, we have $\lb\wt h,x\rb\le a-\la$.
As $\wt\St^*$ is a cone, for any $x\in\wt\St^*$, we have
$\lb\wt h,x\rb\le0$ and $a-\la\ge0$.
Let $h$ be the $d$-dimensional vector that consists of the
first $d$ components of~$\wt h$.
Assume that $h\notin\St$.
Then in the $d$-dimensional plane
$\{x\in\R^{d+1}:x^{d+1}=1\}$ we can select a
$(d-1)$-dimensional plane $L$ that separates
$\wt h$ from $\wt H$.
Consider the $d$-dimensional plane generated by
the origin of $\R^d$ and~$L$, and let $x$ be its normal.
Then $\lb\wt h,x\rb>0$, while
$\sup_{g\in\wt H}\lb g,x\rb\le0$.
Consequently, $x\in \wt H^*$, but then we get a
contradiction with the choice of~$\wt h$.
As a result, $h\in\St$.
Furthermore,
$$
u(W+\lb h,X\rb)
=\inf_{\QQ\in\DDD}\EE_\QQ(W+\lb h,X\rb)
=\inf_{x\in G}\lb\wt h,x\rb
>\la.
$$
As $\la<\la_*$ has been chosen arbitrarily, we conclude that
$\sup_{h\in\St}u(W+\lb h,X\rb)\ge\la_*$.

Let us prove the reverse inequality.
We can assume that $\la_*<\infty$.
Let $x_0\in(\la_*e+\wt\St^*)\cap G$.
Fix $h\in\St$ and set $\wt h=(h,1)$.
Then
$$
u(W+\lb h,X\rb)
=\inf_{x\in G}\lb\wt h,x\rb
\le\lb\wt h,x_0\rb.
$$
We can write $x_0=\la_* e+z_0$ with $z_0\in\wt\St^*$.
Then $\lb\wt h,x_0\rb=\la_*+\lb\wt h,z_0\rb\le\la_*$.
Thus, $\sup_{h\in\St}u(W+\lb h,X\rb)\le\la_*$.
As a result, $u_*=\la_*$.

Let us prove the equality $\argmax_{h\in\St}u(W+\lb h,X\rb)=N$.
In the case $\la_*=\infty$, its left-hand side and its
right-hand side are empty, so it is trivially satisfied.
Assume now that $\la_*<\infty$.
Let $h\in N$.
Using the same arguments as above, we show that $h\in\St$.
For $\wt h=(h,1)$, there exists $a\in\R$ such that,
for any $x\in\la_*e+\wt\St^*$ and any $y\in G$, we have
$\lb h,x\rb\le a\le\lb h,y\rb$.
The same arguments as above show that $a\ge\la_*$.
Consequently,
$$
u(W+\lb h,X\rb)
=\inf_{x\in G}\lb\wt h,x\rb
\ge a
\ge\la_*.
$$

Let $h\in\St$ be such that $u(W+\lb h,X\rb)=\la_*$.
This means that, for $\wt h=(h,1)$, we have
$\inf_{x\in G}\lb\wt h,x\rb\ge\la_*$.
Furthermore, for any $x=\la_*e+z\in\la_*e+\wt\St^*$, we have
$\lb\wt h,x\rb=\lb\wt h,\la_*e\rb+\lb\wt h,z\rb\le\la_*$.
Thus, $\wt h\in\wt N$, which means that $h\in N$.\End

\begin{Example}\rm
\label{GOF2}
{\bf(i)} Let $\St=\R^d$.
Then $\wt\St=\{e\}$,
$\wt\St^*=\{\al e:\al\le0\}$, and
$\la_*=\inf\{x^{d+1}:x\in G_0\}$,
where $G_0=G\cap(\{0\}\times\R)$.
The condition that there exists no $X\in A$ with $u(X)>0$
is equivalent to: $G_0\ne\emp$.
If $G^\circ\cap(\{0\}\times\R)\ne\emp$,
where $G^\circ$ denotes the relative interior of~$G$,
then $N\ne\emp$.
If $G^\circ\cap(\{0\}\times\R)=\emp$, then both
cases $N\ne\emp$ and $N=\emp$ are possible (see Figure~5).

{\bf(ii)} Let $\St=\R_+^d$.
Then $\wt\St=\R_+^d\times\{1\}$, $\wt\St^*=\R_-^{d+1}$,
and $\la_*=\inf\{x^{d+1}:x\in G_-\}$, where
$G_-=G\cap(\R_-^d\times\R)$.\End
\end{Example}

\begin{figure}[h]
\begin{picture}(150,62.5)(-33,-25)
\put(-0.7,-10){\includegraphics{poem.9}}
\put(0,0){\vector(1,0){35}}
\put(0,0){\vector(0,1){35}}
\put(60,0){\vector(1,0){35}}
\put(60,0){\vector(0,1){35}}
\multiput(72,0.5)(0,1.95){15}{\line(0,1){1}}
\multiput(60,29.3)(2,0){6}{\line(1,0){1}}
\put(57.2,23.8){\scalebox{0.9}{$\biggl\{$}}
\put(56,23.7){\small $1$}
\put(0,0){\circle*{1}}
\put(60,0){\circle*{1}}
\put(33,-4){\small $\R^d$}
\put(93,-4){\small $\R^d$}
\put(-4,33){\small $\R$}
\put(56,33){\small $\R$}
\put(13,19){\small $G$}
\put(73,19){\small $G$}
\put(-3,8){\small $e$}
\put(57,8){\small $e$}
\put(72.5,30){\small $\wt h_*$}
\put(71,-4){\small $h_*$}
\put(0.5,-9){\small $\wt\St^*$}
\put(60.5,-9){\small $\wt\St^*$}
\put(7,-22.5){\parbox{80mm}{\small\textbf{Figure~5.}
Existence (right) and nonexistence (left) of an optimal
strategy for the case $\St=\R^d$}}
\end{picture}
\end{figure}

%-------------------------------------------------------
\subsection{Single-Agent Local Optimization}
\label{LOG}

Let $(\Omega,\F,\PP)$ be a probability space,
$u$ be a coherent utility function with the weakly compact
determining set~$\DDD$,
$A\subset L^0$ be a $\DDD$-consistent convex set
containing zero, and $W\in L_s^1(\DDD)$.
The financial interpretation is the same as above.
As opposed to Subsection~\ref{GOG}, we assume that
$\sup_{X\in A,\,\QQ\in\DDD}|\EE_\QQ X|<\infty$
and $A\subseteq L_s^1(\DDD)$.
From the financial point of view, we have a ``big''
investor possessing a capital with the terminal wealth~$W$
and considering several trading opportunities, each of
which is small as compared to~$W$.
Mathematically, we consider the following problem.

\skm
\textbf{Problem (single-agent local optimization):}
Find
$$
u_*=\lim_{\eps\da0}\eps^{-1}
\biggl[\sup_{X\in A}u(W+\eps X)-u(W)\biggr]
$$
and an element $X_*\in A$, for which
\begin{equation}
\label{log1}
\lim_{\eps\da0}\eps^{-1}[u(W+\eps X_*)-u(W)]=u_*.
\end{equation}

\skm
The statement below shows that the problem posed above is
equivalent to the problem of maximizing
$u^c(X;W)=\inf_{\QQ\in\EX_\DDD(Y)}\EE_\QQ X$ over $A$.

\begin{Proposition}
\label{LOG1}
We have $u_*=\sup_{X\in A}u^c(X;W)$.
Furthermore, $X_*$ solves~\eqref{log1} if and only if
$X_*\in\argmax_{X\in A}u^c(X;W)$.
\end{Proposition}

\sksminus
\Proof
The inequality
\begin{align*}
&\limsup_{\eps\da0}\eps^{-1}\biggl(\sup_{X\in A}u(W+\eps X)-u(W)\biggr)\\
&\le\limsup_{\eps\da0}\eps^{-1}\biggl(\sup_{X\in A}
\biggl(\inf_{\QQ\in\EX_\DDD(W)}\EE_\QQ(W+\eps X)-u(W)\biggr)\biggr)\\
&=\limsup_{\eps\da0}\eps^{-1}\sup_{x\in A}\eps
\inf_{\QQ\in\EX_\DDD(W)}\EE_\QQ X\\
&=\sup_{X\in A} u^c(X;W),
\end{align*}
combined with~\cite[Th.~2.16]{C061},
shows that $u_*\le\sup_{X\in A}u^c(X;W)$.
The reverse inequality and the second statement follow
immediately from~\cite[Th.~2.16]{C061}.\End

\skm
Thus, the problem of single-agent local optimization is
equivalent to maximizing another coherent utility over~$A$
(namely, $u^c(\,\cdot\,;W)$).
We will now consider a problem of maximizing a coherent
utility, which we still denote by~$u$, over~$A$ for a
static model with a finite number of assets.
Let $A=\{\lb h,X\rb:h\in\St\}$, where
$X=(X^1,\dots,X^d)\in L_w^1(\DDD)$
and $\St\subset\R^d$ is a convex compact
(here we impose no conditions on the determining set~$\DDD$
of~$u$). Let us introduce the notation (see Figure~6)
\begin{align*}
G&=\cl\{\EE_\QQ X:\QQ\in\DDD\},\\
\St^*&=\{x\in\R^d:\forall h\in\St,\;\lb h,x\rb\le1\},\\
\la_*&=\inf\{\la\ge0:\la\St^*\cap G\ne\emp\},\\
N&=\begin{cases}
\{h:\forall x\in\la_*\St^*,\:\forall y\in G,\;
\lb h,x\rb\le\la_*\le\lb h,y\rb\}&\text{if}\;\;\la_*>0,\\
\emp&\text{if}\;\;\la_*=0.
\end{cases}
\end{align*}
Note that $\la_*<\infty$, $N\ne\emp$, and $N\subseteq\St$.

\begin{figure}[h]
\begin{picture}(150,65)(-61.5,-40)
\put(-21,-21){\includegraphics{poem.10}}
\put(0,0){\circle*{1}}
\put(-1,-4){\small $0$}
\put(35.8,9.5){\small $G$}
\put(3,9){\small $\St$}
\put(-12.5,-11){\small $\St^*$}
\put(-15,19){\small $\la_*\St^*$}
\put(13.5,7){\small $h_*$}
\put(-15,-35){\parbox{67mm}{\small\textbf{Figure~6.}
Solution of the optimization problem.
Here $h_*$ is the optimal~$h$.}}
\end{picture}
\end{figure}

\begin{Theorem}
\label{LOF1}
We have $\sup_{X\in A}u(X)=\la_*$ and
$\argmax_{h\in\St}u(\lb h,X\rb)=N$.
\end{Theorem}

\sksminus
\Proof
Let $\la_*>0$.
For $h\in N$, we have
$$
u(\lb h,X\rb)=\inf_{x\in G}\lb h,x\rb=\la_*.
$$
For $h\in\St\setminus N$, we have
$\sup_{x\in\la_*\St^*}\lb h,x\rb\le\la_*$,
and consequently,
$$
u(\lb h,x\rb)=\inf_{x\in G}\lb h,x\rb<\la_*.
$$
The case $\la_*=0$ is analyzed trivially.\End

%-------------------------------------------------------
\subsection{Liquidity Effects in the NGD Pricing}
\label{UPL}

Let $(\Omega,\F,\PP)$ be a probability space,
$u$ be a coherent utility function with the weakly compact
determining set~$\DDD$,
and $A\subset L^0$ be a convex set containing zero.
We assume that there exists no $X\in A$ with $u(X)>0$.

\begin{Definition}\rm
\label{UPL1}
We define the \textit{upper} and
\textit{lower utility-based NGD price functions} of a
contingent claim~$F$ as
\begin{align*}
\wl V(F,v)&=\sup\{x:\text{the model }(\Omega,\F,\PP,\DDD,A-v(F-x))
\text{ satisfies the NGD}\},\quad v>0,\\
\underline{V}(F,v)&=\inf\{x:\text{the model }(\Omega,\F,\PP,\DDD,A+v(F-x))
\text{ satisfies the NGD}\},\quad v>0.
\end{align*}
From the financial point of view, $v$ means the volume
of a trade.
\end{Definition}

\Remark
If $A$ is a cone, then
\begin{align*}
\wl V(F,\,\cdot\,)&\equiv\inf\{x:\exists X\in A:u(X-F+x)\ge0\},\\
\underline{V}(F,\,\cdot\,)&\equiv\sup\{x:\exists X\in A:u(X+F-x)\ge0\}.
\end{align*}
These are the upper and the lower prices, which were studied
in~\cite[Subsect.~3.6]{C061}. Thus, the investigation
of $\wl V(F,v)$ and $\underline{V}(F,v)$ is meaningful
only if $A$ does not have a cone structure.
This corresponds to the liquidity effects.

\skm
In view of the equality $\underline{V}(F,v)=-\wl V(-F,v)$,
it is sufficient to study only the properties of
$\wl V(F,\,\cdot\,)$.

\begin{Theorem}
\label{UPL2}
Let $F\in L_s^1(\DDD)$.

{\bf(i)} The function $\wl V(F,\,\cdot\,)$ is increasing
and continuous.

{\bf(ii)} We have
$$
\lim_{v\da0}\wl V(F,v)=\sup_{\QQ\in\DDD\cap\RRR}\EE_\QQ F.
$$

{\bf(iii)} We have
$$
\lim_{v\to\infty}\wl V(F,v)\le\sup_{\QQ\in\DDD}\EE_\QQ F.
$$
If $\sup_{X\in A,\,\QQ\in\DDD}|\EE_\QQ X|<\infty$, then
$$
\lim_{v\to\infty}\wl V(F,v)=\sup_{\QQ\in\DDD}\EE_\QQ F.
$$
\end{Theorem}

\Proof
{\bf(i)} It follows from the equality
$$
\sup_{X\in A}u(-v(F-x)+X)=v x+\sup_{X\in A}u(-vF+X)
$$
that $\wl V(F,v)=-v^{-1}f(v)$, where
$f(v)=\sup_{X\in A}u(-vF+X)$.
Note that $f$ is finite due to the NGD and the condition
$F\in L_s^1(\DDD)$.
Fix $v_1,v_2>0$, $\eps>0$, $\al\in[0,1]$
and find $X_1,X_2\in A$ such that
$u(-v_i F+X_i)\ge f(v_i)-\eps$, $i=1,2$. Then
\begin{align*}
f(\al v_1+(1-\al)v_2)
&\ge u(-(\al v_1+(1-\al)v_2)F+\al X_1+(1-\al)X_2)\\
&\ge\al u(-v_1 F+X_1)+(1-\al)u(-v_2 F+X_2)\\
&\ge\al f(v_1)+(1-\al)f(v_2)-\eps.
\end{align*}
Consequently, $f$ is concave.
As $A$ contains zero and the NGD is satisfied, we have $f(0)=0$.
This leads to the desired statement.

{\bf(ii)} By Proposition~\ref{GOG1},
$$
\sup_{X\in\cone A}u(-vF+X)
=\inf_{\QQ\in\DDD\cap\RRR}\EE_\QQ(-v F)
=-v\sup_{\QQ\in\DDD\cap\RRR}\EE_\QQ F,
$$
where ``$\cone$'' denotes the cone hull.
Take $\eps>0$ and find $X_0\in A$, $\al_0\ge0$ such that
$$
u(-F+\al_0 X_0)\ge-\sup_{\QQ\in\DDD\cap\RRR}\EE_\QQ F-\eps.
$$
As the function
$\R_+\ni x\mapsto u(-xF+x\al_0X_0)$ is concave and vanishes at
zero, we have
$$
u(-vF+v\al_0X_0)\ge v\biggl(-\sup_{\QQ\in\DDD\cap\RRR}\EE_\QQ F-\eps\biggr),
\quad v\le1.
$$
As $\eps>0$ has been chosen arbitrarily, we get
$$
\limsup_{v\da0}\wl V(F,v)
=\limsup_{v\da0}(-v^{-1}f(v))
\le\sup_{\QQ\in\DDD\cap\RRR}\EE_\QQ F.
$$
Combining this with the inequality
$$
\sup_{X\in A}u(-vF+X)
\le\sup_{X\in A}\inf_{\QQ\in\DDD\cap\RRR}\EE_\QQ(-vF+X)
=\inf_{\QQ\in\DDD\cap\RRR}\EE_\QQ(-vF)
=-v\sup_{\QQ\in\DDD\cap\RRR}\EE_\QQ F,
$$
we get the desired statement.

{\bf(iii)}
The first statement follows from the inequality
$$
\sup_{X\in A}u(-vF+X)
\ge u(-vF)
=-v\sup_{\QQ\in\DDD}\EE_\QQ F.
$$
The second statement is an obvious consequence of the
equality $\overline{V}(F,v)=-\sup_{X\in A}u(-F+v^{-1}X)$.\End

\skm
\Remarks
\texttt{(i)} If $A$ is a cone, then clearly
$\wl V(F,\,\cdot\,)=const$.

\texttt{(ii)} If $\sup_{X\in A,\,\QQ\in\DDD}|\EE_\QQ X|<\infty$,
then
$$
\wl V(F,\infty)-\underline{V}(F,\infty)
=\sup_{\QQ\in\DDD}F-\inf_{\QQ\in\DDD}F,
$$
which is the length of the NGD price interval in the absence
of a market. The difference
$$
\wl V(F,0)-\underline{V}(F,0)
=\sup_{\QQ\in\DDD\cap\RRR}F-\inf_{\QQ\in\DDD\cap\RRR}F
$$
is the length of the NGD price interval in the presence
of a market. Thus, the ratio
$$
\frac{\wl V(F,0)-\underline{V}(F,0)}%
{\wl V(F,\infty)-\underline{V}(F,\infty)}
$$
measures the ``closeness'' of a new instrument~$F$
to those already existing in the market.

\begin{figure}[!h]
\begin{picture}(150,57.5)(-62,-20)
\put(-2,4.7){\includegraphics{poem.11}}
\put(0,0){\vector(1,0){45}}
\put(0,0){\vector(0,1){35}}
\multiput(-0.5,5)(2,0){9}{\line(1,0){1}}
\multiput(-0.5,29)(2,0){9}{\line(1,0){1}}
\put(-0.5,11.3){\line(1,0){1}}
\put(-0.5,22.7){\line(1,0){1}}
\put(43,-4){\small $\R^d$}
\put(-4,33){\small $\R$}
\put(16.2,15.7){\small $G$}
\put(-13.5,4){\scalebox{0.8}{$\underline{V}(F,\infty)$}}
\put(-12,10){\scalebox{0.8}{$\underline{V}(F,0)$}}
\put(-12,22){\scalebox{0.8}{$\wl V(F,0)$}}
\put(-13.5,28){\scalebox{0.8}{$\wl V(F,\infty)$}}
\put(-9.5,-15){\parbox{55mm}{\small\textbf{Figure~7.}
The form of $\wl V(F,0)$, $\wl V(F,\infty)$,
$\underline{V}(F,0)$, and $\underline{V}(F,\infty)$}}
\end{picture}
\end{figure}

\begin{Example}\rm
\label{UPL3}
Consider a static model with a finite number of assets,
i.e. $A=\{\lb h,X\rb:h\in\St\}$, where
$X=(X^1,\dots,X^d)\in L_w^1(\DDD)$ and
$\St\subset\R^d$ is a convex bounded set.
Assume that $\St$ contains a neighborhood of zero.
Consider the generator $G=\{\EE_\QQ(X,F):\QQ\in\DDD\}$.
Then
\begin{align*}
\wl V(F,0)&=\sup\{x^{d+1}:x^1=\dots=x^d=0,\;x\in G\},\\
\wl V(F,\infty)&=\sup\{x^{d+1}:x\in G\}.
\end{align*}
Note that these values do not depend on~$\St$!\End
\end{Example}

%========================================================
\section{Optimality Pricing}
\label{OP}

%--------------------------------------------------------
\subsection{Agent-Independent Optimality Pricing}
\label{AIOP}

Consider the model of~\cite[Subsect.~3.2]{C061}.
Thus, we are given a probability space $(\Omega,\F,\PP)$,
a convex weakly compact set $\RD\subset\PPP$,
an $L^1$-closed convex set $\PD\subset\RD$,
and a convex set $A\subset L^0$.
Assume that $R_*<\infty$, where $R_*=\sup_{X\in A}\RAROC(X)$.
It follows from~\cite[Th.~3.10]{C061} that
$$
R_*=\inf\biggl\{R\ge0:\Bigl(\frac{1}{1+R}\,\PD
+\frac{R}{1+R}\,\RD\Bigr)\cap\RRR\ne\emp\biggr\}
$$
and $\DDD_*\cap\RRR\ne\emp$, where
$$
\DDD_*=\frac{1}{1+R_*}\,\PD+\frac{R_*}{1+R_*}\,\RD.
$$

\begin{figure}[h]
\begin{picture}(150,62.5)(-45,-10)
\put(0,0){\includegraphics{poem.12}}
\put(27,24){\small $\PD$}
\put(28,32){\small $\DDD_*$}
\put(27,40){\small $\RD$}
\put(63.5,34){\small $\RRR$}
\put(48.5,24.5){\small $\DDD_*\!\!\cap\!\RRR$}
\put(5,-7.5){\small\textbf{Figure~8.}
The structure of $\DDD_*\cap\RRR$}
\end{picture}
\end{figure}

\begin{Definition}\rm
\label{AIOP1}
An \textit{agent-independent NBC price} of
a contingent claim~$F$ is a real number $x$ such that
$$
\sup_{X\in A+A(x)}\RAROC(X)=\sup_{X\in A}\RAROC(X),
$$
where $A(x)=\{h(F-x):h\in\R\}$.

The set of the NBC prices will be denoted by $\INBC(F)$.
\end{Definition}

This pricing technique corresponds to the agent-independent
optimization. A price~$x$ is fair if adding to the market
a new instrument with the initial price~$x$ and the
terminal price~$F$ does not increase the optimal value in
the optimization problem.

\begin{Proposition}
\label{AIOP2}
For $F\in L_s^1(\RD)$,
$$
\INBC(F)=\{\EE_\QQ F:\QQ\in\DDD_*\cap\RRR\}.
$$
\end{Proposition}

\sksmin
\Proof
If $x\in\INBC(F)$, then, by~\cite[Th.~3.10]{C061},
there exists $\QQ\in\DDD_*\cap\RRR(A+A(x))$.
This means that $\QQ\in\DDD_*\cap\RRR$ and $\EE_\QQ F=x$.

Conversely, if $x=\EE_\QQ F$ with some $\QQ\in\DDD_*\cap\RRR$,
then, for any $X+h(F-x)\in A+A(x)$, we have
$\EE_\QQ X\le0$, so that $\QQ\in\RRR(A+A(x))$.
Due to~\cite[Th.~3.10]{C061},
$\sup_{X\in A+A(x)}\RAROC(X)\le R_*$.\End

\skm
The following statement yields a more definite representation
of $\DDD_*\cap\RRR$.

\begin{Proposition}
\label{AIOP3}
If $X_*\in\argmax_{X\in A}\RAROC(X)$, then
\begin{equation}
\label{aiop1}
\DDD_*\cap\RRR
=\biggl(\frac{1}{1+R_*}\,\EX_\PD(X_*)+\frac{R_*}{1+R_*}\,\EX_\RD(X_*)\biggr)\cap\RRR.
\end{equation}
\end{Proposition}

\sksmin
\Proof
Take
$$
\QQ=\frac{1}{1+R_*}\,\QQ_1+\frac{R_*}{1+R_*}\,\QQ_2\in\DDD_*\cap\RRR.
$$
We have
$$
\inf_{\QQ\in\PD}\EE_\QQ X_*+R_*\inf_{\QQ\in\RD}\EE_\QQ X_*
\le\EE_{\QQ_1}X_*+R_*\EE_{\QQ_2}X_*
\le0
$$
(the second inequality follows from the inclusion
$\QQ\in\RRR$).
Combining this with the equality
$$
\RAROC(X_*)
=\frac{\inf_{\QQ\in\PD}\EE_\QQ X_*}{-\inf_{\QQ\in\RD}\EE_\QQ X_*}
=R_*,
$$
we get
$$
\inf_{\QQ\in\PD}\EE_\QQ X_*+R_*\inf_{\QQ\in\RD}\EE_\QQ X_*
\le\EE_{\QQ_1}X_*+R_*\EE_{\QQ_2}X_*.
$$
This means that $\QQ_1\in\EX_\PD(X_*)$ and
$\QQ_2\in\EX_\RD(X_*)$.\End

\skm
As a corollary, if $\EX_\PD(X_*)$ and $\EX_\RD(X_*)$ are
singletons (this is true, for instance, if $\PD=\{\PP\}$,
$\RD$ is the determining set of Weighted V@R, and
$X_*$ has a continuous distribution; see~\cite{C05e}),
then $\RRR$ can be removed from~\eqref{aiop1}, i.e.
$$
\DDD_*\cap\RRR
=\frac{1}{1+R_*}\,\EX_\PD(X_*)+\frac{R_*}{1+R_*}\,\EX_\RD(X_*).
$$
But in general this equality might be violated as shown
by the example below.

\begin{Example}\rm
\label{AIOP4}
Let $\PD=\{\PP\}$, $\RD$ be the determining set of
Tail V@R with $\la<1/2$,
and $X^1$, $X^2$ be independent random variables
with $\PP(X^1=-1)=\PP(X^1=2)=1/2$, $\PP(X^2=\pm1)=1/2$.
Let $A=\{h^1X^1+h^2X^2:h^i\in\R\}$.
For any $(h^1,h^2)$ with $h^1\ge0$, we have
$$
\inf_{\QQ\in\RD}\EE_\QQ(h^1X^1+h^2X^2)
\le\EE_\PP Z(h^1X^1+h^2X^2)
=h^1\EE_\PP ZX^1
=\inf_{\QQ\in\RD}\EE_\QQ X^1,
$$
where $Z=2I(X=-1)$.
Combining this with the equality
$\EE_\PP(h^1X^1+h^2X^2)=h^1\EE_\PP X^1$,
we get that $X^1\in\argmax_{X\in A}\RAROC(X)$.
On the other hand, there exists $\QQ\in\EX_\RD(X^1)$,
for which $\EE_\QQ X^2\ne0$. Thus, the set
$\frac{1}{1+R_*}\,\EX_\PD(X^1)
+\frac{R_*}{1+R_*}\,\EX_\RD(X^1)$
contains measures that do not belong to~$\RRR$.\End
\end{Example}

\Remark
One of techniques for pricing in incomplete markets
consists in finding the representative of the set of
risk-neutral measures that is the closest one to~$\PP$
in some sense (typically the relative entropy or some
other measure of distance is minimized).
Note that the set $\DDD_*\cap\RRR$ is exactly the set of
measures~$\QQ$ from~$\RRR$ that are the closest ones
to~$\PD$, the ``distance'' being measured by
$$
\inf\biggl\{R:\exists\QQ_1\in\PD,\,\QQ_2\in\RD:
\frac{1}{1+R}\,\QQ_1+\frac{R}{1+R}\,\QQ_2=\QQ\biggr\}.
$$

\skm
We will now study the problem for a static model with
a finite number of assets.
Let $A=\{\lb h,S_1-S_0\rb:h\in\St\}$, where
$S_0\in\R^d$, $S_1^1,\dots,S_1^d\in L_w^1(\RD)$, and
$\St\subseteq\R^d$ is a closed convex cone.
Assume that $0<R_*<\infty$.
Let $F\in L_w^1(\DDD)$ be a contingent claim.
Let us introduce the notation (see Figure~9)
\begin{align*}
\St^*&=\{x\in\R^d:\forall h\in\St,\:\lb h,x\rb\ge0\},\\
\wt\St^*&=\St^*\times\{0\},\\
\wt E&=\cl\{\EE_\QQ(S_1,F):\QQ\in\PD\},\\
\wt G&=\cl\{\EE_\QQ(S_1,F):\QQ\in\RD\},\\
\wt D&=\wt G+\wt\St^*,\\
\wt D_R&=\frac{1}{1+R}\,\wt E+\frac{R}{1+R}\,\wt D.
\end{align*}

\begin{figure}[h]
\begin{picture}(150,85)(-75,-53)
\put(-40,-24){\includegraphics{poem.13}}
\put(-1.7,-1.5){\small $\wt E$}
\put(8.3,19.5){\small $\wt G$}
\put(30.3,19.7){\small $\wt D$}
\put(29.5,13.7){\small $\wt D_R$}
\put(28.5,7.8){\small $\wt D_{R_*}$}
\put(48,-34){\small $\R^d$}
\put(-39,28){\small $\R$}
\put(-14,-34){\small $S_0$}
\put(-55,-6){\small $\INBC(F)$}
\put(-35,-30){\vector(1,0){85}}
\put(-35,-30){\vector(0,1){60}}
\multiput(-11.7,-29.5)(0,2){25}{\line(0,1){1}}
\multiput(-34.5,-5)(2,0){12}{\line(1,0){1}}
\put(-35,-5){\circle*{1}}
\put(-24,-48){\parbox{59mm}{\small\textbf{Figure~9.}
The form of $\INBC$. Here $\INBC(F)$ consists of
one point.}}
\end{picture}
\end{figure}

\begin{Theorem}
\label{AIOPF1}
We have
\begin{align}
\label{saopf1}
R_*&=\inf\{R>0:\wt D_R\cap(\{S_0\}\times\R)\ne\emp\},\\
\label{saopf2}
\INBC(F)&=\{x:(S_0,x)\in\wt D_{R_*}\}.
\end{align}
\end{Theorem}

\sksmin
\Proof
Denote
\begin{align*}
E&=\cl\{\EE_\QQ S_1:\QQ\in\PD\},\\
G&=\cl\{\EE_\QQ S_1:\QQ\in\RD\},\\
D&=G+\St^*,\\
D_R&=\frac{1}{1+R}\,E+\frac{R}{1+R}\,D.
\end{align*}
Note that
$E=\pr_{\R^d}\wt E$,
$G=\pr_{\R^d}\wt G$,
$\St^*=\pr_{\R^d}\wt\St^*$,
and consequently,
$D=\pr_{\R^d}\wt D$,
$D_R=\pr_{\R^d}\wt D_R$.
Combining this with the results of Subsection~\ref{OG},
we get
$$
R_*
=\inf\{R>0:D_R\ni S_0\}
=\inf\{R>0:\wt D_R\cap(\{S_0\}\times\R)\ne\emp\}.
$$
Furthermore, for any $x\in\R$,
$$
\sup_{A+A(x)}\RAROC(X)
=\inf\{R>0:\wt D_R\ni(S_0,x)\}.
$$
This, combined with~\eqref{saopf1}, proves~\eqref{saopf2}.\End

\skm
To conclude this subsection, we find the form of $\INBC(F)$
in the Gaussian case.

\begin{Example}\rm
\label{AIOPF2}
Consider the setting of~\cite[Ex.~3.13]{C061}.
Clearly, $R_*$ is the solution of the equation
$\lb S_0-a,C^{-1}(S_0-a)\rb=\frac{\ga^2R_*^2}{1+R_*^2}$
(cf. Example~\ref{OF3}).
This, combined with the form of $\INGDR(F)$
found in~\cite[Ex.~3.13]{C061}, shows that
$\INBC(F)$ consists of a unique point $\lb b,S_0-a\rb+\EE F$.
Let us remark that this value coincides with the fair price
of $F$ obtained as a result of the mean-variance hedging.
Note that this value does not depend on~$u$!\End
\end{Example}

%--------------------------------------------------------
\subsection{Single-Agent Optimality Pricing}
\label{SAOP}

Let $(\Omega,\F,\PP)$ be a probability space,
$u$ be a coherent utility function with the weakly compact
determining set~$\DDD$,
$A\subseteq L^0$ be a $\DDD$-consistent convex set
containing zero, and $W\in L_s^1(\DDD)$.
The financial interpretation is the same as in
Subsection~\ref{GOG}.

\begin{Definition}\rm
\label{SAOP1}
A \textit{single-agent NBC price} of a contingent claim~$F$
is a real number~$x$ such that
$$
\max_{X\in A,\,h\in\R}u(W+X+h(F-x))=u(W).
$$

The set of the NBC prices will be denoted by $\INBC(F)$.
\end{Definition}

This pricing technique corresponds to the global
single-agent optimization. A price~$x$ is fair if adding
to the market a new instrument with the initial price~$x$
and the terminal price~$F$ does not increase the optimal
value in the optimization problem.

\begin{Theorem}
\label{SAOP2}
For $F\in L_s^1(\DDD)$,
\begin{equation}
\label{saop1}
\INBC(F)=\{\EE_\QQ F:\QQ\in\EX_\DDD(W)\cap\RRR\}.
\end{equation}
\end{Theorem}

\Remark
The set of the NBC prices is nonempty only if $W$ is optimal
in the sense that $\max_{X\in A}u(W+X)=u(W)$.
However, if $W$ is not optimal, then, as seen from the
proof of Theorem~\ref{SAOP2}, $\EX_\DDD(W)\cap\RRR=\emp$,
so that~\eqref{saop1} still holds.

\skm
\texttt{Proof of Theorem~\ref{SAOP2}.}
As $A$ contains zero and
the function $\R_+\ni\al\mapsto u(W+\al X)$ is concave
for a fixed~$X$,
the condition $x\in\INBC(F)$ is equivalent to:
$$
\max_{X\in\cone A,\,h\in\R}u(W+X+h(F-x))=u(W).
$$
By Proposition~\ref{GOG1}, this is equivalent to:
$$
\inf_{\QQ\in\DDD\cap\RRR(A+A(x))}\EE_\QQ W
=\inf_{\QQ\in\DDD}\EE_\QQ W,
$$
where $A(x)=\{h(F-x):h\in\R\}$.
Clearly, the latter condition is equivalent to:
$\EX_\DDD(W)\cap\RRR(A+A(x))\ne\emp$.
It is easy to verify that this is equivalent to:
$x=\EE_\QQ F$ for some $\QQ\in\EX_\DDD(W)\cap\RRR$.\End

\skm
Let us now provide a geometric representation of $\INBC(F)$
(see Figure~10).
Assume that $u(W)=\max_{X\in A}u(W+X)$
(the reasoning used above shows that this is equivalent
to: $\EX_\DDD(W)\cap\RRR\ne\emp$).
Consider the generator $G=\{\EE_\QQ(F,W):\QQ\in\DDD\cap\RRR\}$
and the function $f(x)=\inf\{y:(x,y)\in G\}$
(we set $\inf\emp=+\infty$).

\begin{figure}[h]
\begin{picture}(150,92.5)(-45.5,-50)
\put(11.7,-23){\includegraphics{poem.14}}
\put(0,0){\vector(1,0){70}}
\put(0,0){\vector(0,1){40}}
\multiput(12,-0.5)(0,2){11}{\line(0,1){1}}
\multiput(20,-0.5)(0,2){9}{\line(0,1){1}}
\multiput(35,0.5)(0,2){3}{\line(0,1){1}}
\multiput(50,-0.5)(0,2){9}{\line(0,1){1}}
\multiput(58,-0.5)(0,2){11}{\line(0,1){1}}
\put(35,0){\circle*{1}}
\put(68,-3){\small $x$}
\put(-3,38.5){\small $y$}
\put(33,14.5){\small $G$}
\put(33,28.5){\small $D$}
\put(23,11){\scalebox{0.8}{$f(x)$}}
\put(24,-6){\scalebox{0.8}{$I_{\text{\sl NBC}(\DDD,A,W)}(F)$}}
\put(23.5,-10){\scalebox{0.8}{$=\!I_{\text{\sl NBC}(\DDD,W)}(F)$}}
\put(25.5,-17){\scalebox{0.8}{$I_{\text{\sl NGD}(\DDD,A)}(F)$}}
\put(27,-24){\scalebox{0.8}{\small $I_{\text{\sl NGD}(\DDD)}(F)$}}
\put(-2.5,-40){\parbox{74mm}{\small\textbf{Figure~10.}
Comparison of various price intervals.
Here $G=\{\EE_\QQ(F,W):\QQ\in\DDD\cap\RRR\}$
and $D=\{\EE_\QQ(F,W):\QQ\in\DDD\}$.
In this example,
$I_{\text{\sl NBC}(\DDD,A,W)}(F)=I_{\text{\sl NBC}(\DDD,W)}(F)$.}}
\end{picture}
\end{figure}

\begin{Corollary}
\label{SAOP3}
For $F\in L_s^1(\DDD)$,
$$
\INBC(F)=\argmin_{x\in\R}f(x).
$$
\end{Corollary}

\sksmin
\Proof
It is sufficient to note that
$$
\min_{x\in\R}f(x)=\min_{\QQ\in\DDD\cap\RRR}\EE_\QQ W=u(W)
$$
and
$$
f(x)=\inf\{\EE_\QQ W:\QQ\in\DDD\cap\RRR(A+A(x))\},
$$
where $A(x)=\{h(F-x):h\in\R\}$.
Thus, $x\in\argmin_{x\in\R}f(x)$ if and only if
$\EX_\DDD(W)\cap\RRR(A+A(x))\ne\emp$, which, in view
of Theorem~\ref{SAOP2}, is equivalent to the inclusion
$x\in\INBC(F)$.\End

\skm
Assume that $W$ is optimal in the sense that
\begin{equation}
\label{saop2}
u(W)=\max_{X\in A}u(W+X)
\end{equation}
and suppose moreover that the set
$I_{\text{\sl NBC}(\DDD,W)}(F)$ of the NBC prices based on
$\DDD$ and $W$ (with $A=0$) consists of one point~$x_0$
(this condition is satisfied if the set
$D=\{\EE_\QQ(F,W):\QQ\in\DDD\}$ is strictly convex;
see Figure~10).
It is seen from the proof of Theorem~\ref{SAOP2} that
condition~\eqref{saop2} is equivalent to:
$\EX_\DDD(W)\cap\RRR\ne\emp$.
Then it follows from Theorem~\ref{SAOP2} that
$I_{\text{\sl NBC}(\DDD,A,W)}(F)\ne\emp$
(we assume that $F\in L_s^1(\DDD)$).
Clearly, $I_{\text{\sl NBC}(\DDD,A,W)}(F)\subseteq
I_{\text{\sl NBC}(\DDD,W)}(F)$.
As a result, $I_{\text{\sl NBC}(\DDD,A,W)}(F)=\{x_0\}$.
So, in this situation $A$ can be eliminated.
This situation occurs naturally as shown, in particular,
by the example below.

\begin{Example}\rm
\label{SAOP4}
Let $u$ be a law invariant coherent utility function that
is finite on Gaussian random variables.
Assume that $u(W)=\max_{X\in A}u(W+X)$ and
that $(W,F)$ has a Gaussian distribution.

There exists $\ga>0$ such that, for a Gaussian random
variable~$\xi$ with mean~$m$ and variance~$\si^2$, we have
$u(\xi)=m-\ga\si$.
Clearly, $\INBC(F)\subseteq J$, where $J$ is the NBC
price based on $\DDD$ and $W$ with $A=0$.
Using Corollary~\ref{SAOP3}, we deduce that
$J$ consists of a single point
$\EE F-\ga\frac{\cov(F,W)}{(\var W)^{1/2}}$.
As $\INBC(F)$ is nonempty, it consists of the same point.\End
\end{Example}

%--------------------------------------------------------
\subsection{Multi-Agent Optimality Pricing}
\label{MAOP}

Let $(\Omega,\F,\PP)$ be a probability space,
$u_1,\dots,u_N$ be coherent utility functions with the
weakly compact determining sets $\DDD_1,\dots,\DDD_N$,
$A\subseteq L^0$ be a convex set containing zero,
and $W_1\in L_s^1(\DDD_1),\dots,W_N\in L_s^1(\DDD_N)$.
From the financial point of view,
$u_n$, $A$, and $W_n$ are the coherent utility function,
the set of attainable P\&Ls,
and the terminal endowment of the $n$-th agent, respectively.
We will assume that there exists a set
$A'\subseteq\bigcap_n L_s^1(\DDD_n)\cap A$ such that,
for any $n$, $\DDD_n\cap\RRR=\DDD_n\cap\RRR(A')$.
We also assume that each $W_n$ is optimal in the sense
that $u_n(W_n)=\max_{X\in A}u_n(W_n+X)$.

\begin{Definition}\rm
\label{MAOP1}
A real number $x$ is a \textit{multi-agent NBC price}
of a contingent claim~$F$ if there exists no element
$X\in A+\{h(F-x):h\in\R\}$ such that
$u_n(W_n+X)>u_n(W_n)$ for any~$n$.

The set of the NBC prices will be denoted by $\INBC(F)$.
\end{Definition}

From the financial point of view, a price~$x$ is fair if
adding to the market a new instrument with the initial price~$x$
and the terminal price~$F$ does not produce a trading
opportunity that is attractive to all the agents.

\begin{Theorem}
\label{MAOP2}
For $F\in\bigcap_n L_s^1(\DDD_n)$,
$$
\INBC(F)
=\conv_{n=1}^N I_{\text{\sl NBC}(\DDD_n,A,W_n)}(F)
=\bigl\{\EE_\QQ F:\QQ\in\conv_{n=1}^N(\EX_{\DDD_n}(W_n)\cap\RRR)\bigr\},
$$
where $I_{\text{\sl NBC}(\DDD_n,A,W_n)}(F)$ is the interval of
the single-agent NBC prices based on $\DDD_n$, $A$,~$W_n$.
\end{Theorem}

\sksminus
\Proof
Let $x\in\INBC(F)$.
Fix $X_1,\dots,X_M\in A'$.
It follows from the weak continuity of the
maps $\DDD_n\ni\QQ\mapsto\EE_\QQ(X_1,\dots,X_M,F)$ that,
for each~$n=1,\dots,N$, the set
$G_n=\{\EE_\QQ(X_1,\dots,X_M,F-x):\QQ\in\EX_n\}$,
where $\EX_n=\EX_{\DDD_n}(W_n)$, is compact.
Clearly, $G_n$ is convex.
Suppose that
$$
\bigl(\conv_{n=1}^N G_n\bigr)\cap
(\R_-^M\times\{0\})=\emp.
$$
Then there exists $h\in\R^{M+1}$ such that
$h_1,\dots,h_M\ge0$ and
$\inf_{x\in G_n}\lb h,x\rb>0$ for each~$n$.
This means that $\inf_{\QQ\in\EX_n}\EE_\QQ Y>0$
for each~$n$,
where $Y=h^1X_1+\dots+h^mX_m+h^{m+1}(F-x)$.
Employing~\cite[Th.~2.16]{C061}, we conclude that there
exists $\eps>0$ such that $u(W_n+\eps Y)>u(W_n)$ for any~$n$.

The obtained contradiction shows that,
for any $X_1,\dots,X_M\in A'$, the set
\begin{align*}
B(X_1,\dots,X_M)
&=\biggl\{\al_1,\dots,\al_N,\QQ_1,\dots,\QQ_N
\in S\times\prod_{n=1}^N\EX_n:
\sum_{n=1}^N\al_n\EE_{\QQ_n}F=x\\
&\hspace*{7mm}\text{and }\forall n=1,\dots,N,\:\forall m=1,\dots,M,\;\EE_{\QQ_n}X_m\le0\biggr\},
\end{align*}
where $S=\bigl\{\al_1,\dots,\al_N\ge0:\sum_{n=1}^N\al_n=1\bigr\}$,
is nonempty.
As the map $\EX_n\ni\QQ\mapsto\EE_\QQ X$ is
weakly continuous for each $X\in L_s^1(\DDD_n)$,
the set $B(X_1,\dots,X_M)$ is closed with respect to
the product of weak topologies.
Furthermore, any finite intersection of sets of this
form is nonempty. Tikhonov's theorem
ensures that $S\times\prod_n\EX_n$ is compact.
Consequently, there exists a collection
$\al_1,\dots,\al_N,\QQ_1,\dots,\QQ_N$
that belongs to each~$B$ of this form.
Then $\EE_{\QQ_n}X\le0$ for any $n$ and any $X\in A'$,
which means that $\QQ_n\in\EX_n\cap\RRR$.
Thus, the measure $\QQ=\sum_n\al_n\QQ_n$ belongs to
$\conv_n(\EX_n\cap\RRR)$ and $\EE_\QQ F=x$.

Now, let $x=\EE_\QQ F$ with $\QQ=\sum_n\al_n\QQ_n$,
$\QQ_n\in\EX_n\cap\RRR$.
Suppose that there exist $X\in A$, $h\in\R$ such that, for
$Y=X+h(F-x)$, we have $u_n(W_n+Y)>u_n(W_n)$ for each~$n$.
Due to the concavity of the function
$\al\mapsto u_n(W_n+\al Y)$, we get
\begin{align*}
u_n(W_n+Y)-u_n(W_n)
&\le\limsup_{\eps\da0}\eps^{-1}(u_n(W_n+\eps Y)-u_n(W_n))\\
&\le\limsup_{\eps\da0}\eps^{-1}\biggl(
\inf_{\QQ\in\EX_n}\EE_\QQ(W_n+\eps Y)-u_n(W_n)\biggr)
=\inf_{\QQ\in\EX_n}\EE_\QQ Y.
\end{align*}
Consequently, $\EE_{\QQ_n}Y>0$ for each~$n$,
and therefore, $\EE_\QQ Y>0$.
But, on the other hand, $\QQ\in\RRR$, and therefore,
$\EE_\QQ Y\le\EE_\QQ h(F-x)=0$.
The contradiction shows that $x\in\INBC(F)$.\End

%========================================================
\section{Equilibrium}
\label{E}

%--------------------------------------------------------
\subsection{Unconstrained Equilibrium}
\label{EC}

Let $(\Omega,\F,\PP)$ be a probability space,
$u_1,\dots,u_N$ be coherent utility functions with the
weakly compact determining sets $\DDD_1,\dots,\DDD_N$,
$A_1,\dots,A_N\subseteq L^0$ be convex cones
such that $A_n$ is $\DDD_n$-consistent for each~$n$,
and let $W_1\in L_s^1(\DDD_1),\dots,W_N\in L_s^1(\DDD_N)$.
From the financial point of view,
$u_n$, $A_n$, and $W_n$ are the coherent utility function,
the ``personal'' set of attainable P\&Ls,
and the terminal endowment of the $n$-th agent, respectively.
Let us introduce the notation $\DDD=\bigcap_n\DDD_n$.

\begin{Definition}\rm
\label{EC1}
The \textit{maximal overall utility} is defined as
$$
M=\sup_{\vbox{
\vspace{1mm}
\hbox{\hspace*{10mm}$\scriptstyle X_n\in A_n$}
\vspace{-1mm}
\hbox{$\scriptstyle Y_n\in L_s^1(\DDD)\,:\,\sum_n\!Y_n=0$}}}
\hspace{1mm}\sum_{n=1}^N u_n(W_n+X_n+Y_n),
$$
where the sum is understood as $-\infty$ if any of the
summands equals $-\infty$.
\end{Definition}

\begin{Proposition}
\label{EC2}
We have
$$
M=\inf_{\QQ\in\bigcap_n\!\DDD_n\cap\RRR(A_n)}\EE_\QQ W,
$$
where $W:=\sum_n W_n$ and $\inf\emp:=\infty$.
\end{Proposition}

\begin{Lemma}
\label{EC3}
Let $u_1,\dots,u_N$ be coherent utility functions with
the weakly compact determining sets $\DDD_1,\dots,\DDD_N$.
Then, for any $X\in L^\infty$,
\begin{equation}
\label{ec1}
\sup_{X_n\in L^\infty\,:\,\sum_n\!X_n=X}\hspace{1mm}
\sum_{n=1}^N u_n(X_n)
=\inf_{\QQ\in\bigcap_n\!\DDD_n}\EE_\QQ X.
\end{equation}
\end{Lemma}

\Remark
The left-hand side of~\eqref{ec1} is called the
\textit{convex convolution} or the \textit{sup-convolution}
of $u_1,\dots,u_N$ (see~\cite{BE05},
\cite[Sect.~5.2]{D05}).
Thus, Lemma~\ref{EC3} states that it is a coherent utility
function with the determining set $\bigcap_n\DDD_n$ if
$\bigcap_n\DDD_n\ne\emp$ and it is identically equal to
$+\infty$ if $\bigcap_n\DDD_n=\emp$.

\skm
\texttt{Proof of Lemma~\ref{EC3}. }
In the case, where $\bigcap_n\DDD_n\ne\emp$,
this statement follows by induction from a result
proved in~\cite[Sect.~5.2]{D05}.

Assume now that $\bigcap_n\DDD_n=\emp$.
Find $m$ such that $\bigcap_{n=1}^m\DDD_n\ne\emp$, while
$\bigcap_{n=1}^{m+1}\DDD_n=\emp$.
By the Hahn-Banach theorem, there exists $Z\in L^\infty$
such that
$$
\sup_{\QQ\in\DDD_{m+1}}\EE_\QQ Z
<0
<\inf_{\QQ\in\bigcap_{n=1}^m\DDD_n}\EE_\QQ Z.
$$
According to the part of the lemma that has already been
proved, there exist $Z_1,\dots,Z_m\in L^\infty$ such that
$\sum_{n=1}^m Z_n=Z$ and $\sum_{n=1}^m u_n(Z_n)>0$.
Then $u_1(Z_1)+\dots+u_m(Z_m)+u_{m+1}(-Z)>0$.
Consequently, the left-hand side of~\eqref{ec1} is
identically equal to~$\infty$.\End

\skm
\texttt{Proof of Proposition~\ref{EC2}.}
For any $X_1\in A_1,\dots,X_N\in A_N,
Y_1,\dots,Y_N\in L_s^1(\DDD)$ such that $\sum_n Y_n=0$,
and any $\QQ\in\bigcap_n\DDD_n\cap\RRR(A_n)$, we have
$$
\sum_{n=1}^N u_n(W_n+X_n+Y_n)
\le\sum_{n=1}^N\EE_\QQ(W_n+X_n+Y_n)
\le\sum_{n=1}^N\EE_\QQ(W_n+Y_n)
=\EE_\QQ W
$$
(to get the second inequality, we used the inclusions
$W_n\in L^1(\QQ)$, $Y_n\in L^1(\QQ)$).
Consequently,
$$
M\le\inf_{\QQ\in\bigcap_n\!\DDD_n\cap\RRR(A_n)}
\EE_\QQ W.
$$

Let us prove the reverse inequality.
Clearly, it is sufficient to prove it for bounded $W_n$
(since arbitrary $W_n\in L_s^1(\DDD_n)$ can be approximated
by bounded ones).
Proposition~\ref{GOG1} and Lemma~\ref{EC3} combined together
yield
\begin{align*}
M&\ge\sup_{\vbox{
\vspace{1mm}
\hbox{\hspace*{7mm}$\scriptstyle X_n\in A_n$}
\vspace{-1mm}
\hbox{$\scriptstyle Y_n\in L^\infty\,:\,\sum_n\!Y_n=0$}}}
\hspace{1mm}\sum_{n=1}^N u_n(W_n+X_n+Y_n)\\
&=\sup_{Y_n\in L^\infty\,:\,\sum_n\!Y_n=0}\hspace{1mm}
\sum_{n=1}^N\inf_{\QQ\in\DDD_n\cap\RRR(A_n)}\EE_\QQ(W_n+Y_n)\\
&=\inf_{\QQ\in\bigcap_n\!\DDD_n\cap\RRR(A_n)}
\EE_\QQ W.
\end{align*}

\skm
The following example shows that the restriction
$Y_n\in L_s^1(\DDD)$ in the definition of~$M$
is essential for Proposition~\ref{EC2} and cannot be
eliminated.

\begin{Example}\rm
\label{EC4}
Let $N=2$, $\DDD_1=\DDD_2=\{\PP\}$,
$A_n=\{a\xi_n+Y:a\in\R,\,Y\in L^\infty,$ $\EE_\PP Y=0\}$,
where $\EE_\PP\xi_n^+=\EE_\PP\xi_n^-=\infty$, $n=1,2$,
$\xi_1+\xi_2=1$, and $W_1=W_2=0$.
Take $X_n=\xi_n$, $Y_n=1/2-\xi_n$, $n=1,2$. Then
$Y_1+Y_2=0$ and $W_n+X_n+Y_n=1/2$, so that
$$
\sup_{\vbox{
\vspace{1mm}
\hbox{\hspace*{7mm}$\scriptstyle X_n\in A_n$}
\vspace{-1mm}
\hbox{$\scriptstyle Y_n\in L^0\,:\,\sum_n\!Y_n=0$}}}
\hspace{1mm}\sum_{n=1}^N u_n(W_n+X_n+Y_n)
=\infty.
$$
On the other hand, $\bigcap_n\DDD_n\cap\RRR(A_n)=\{\PP\}$
and $\EE_\PP W=0$.\End
\end{Example}

We now pass on to various definitions of equilibrium.
From the financial point of view, the Pareto-type
equilibrium corresponds to the global optimum, while
the Arrow-Debreu-type equilibrium corresponds to the
competitive optimum.

\begin{Definition}\rm
\label{EC5}
A \textit{Pareto-type equilibrium} is a collection
$(X,Y)=$ $(X_1,\dots,X_N,Y_1,\dots,Y_N)$ such that
\begin{mitemize}
\item[(a)] $X_n\in A_n$;
\item[(b)] $Y_n\in L_s^1(\DDD)$, $\sum_n Y_n=0$;
\item[(c)] there do not exist $(X',Y')$ satisfying (a), (b)
and such that
\begin{align*}
\forall n,\;u_n(W_n+X'_n+Y'_n)\ge u_n(W_n+X_n+Y_n),\\
\exists n:u_n(W_n+X'_n+Y'_n)>u_n(W_n+X_n+Y_n).
\end{align*}
\end{mitemize}
\end{Definition}

It is easy to see from the translation invariance
property ($u_n(X+m)=u_n(X)+m$) that condition~(c) is
equivalent to:
\begin{mitemize}
\item[(c')] $\sum_n u_n(W_n+X_n+Y_n)=M$.
\end{mitemize}

\begin{Definition}\rm
\label{EC6}
An \textit{Arrow-Debreu-type equilibrium}
is a collection $(X,Y,\QQ)$, where $\QQ\in\PPP$, such that
\begin{mitemize}
\item[(a)] $X_n\in A_n$;
\item[(b)] $Y_n\in L_s^1(\DDD)$, $\sum_n Y_n=0$,
$\EE_\QQ Y_n=0$ (so that automatically $Y_n\in L^1(\QQ)$);
\item[(c)] for any $n$,
$$
u_n(W_n+X_n+Y_n)=\max_{\vbox{
\vspace{1mm}
\hbox{\hspace*{8mm}$\scriptstyle\xi\in A_n$}
\vspace{-1mm}
\hbox{$\scriptstyle\eta\in L^1(\QQ)\,:\,\EE_\QQ\eta=0$}}}
\hspace{1mm}u_n(W_n+\xi+\eta).
$$
\end{mitemize}
\end{Definition}

Below the notation $Y'\sim Y$ for $N$-dimensional
random vectors $(Y'_1,\dots,Y'_N)$ and $(Y_1,\dots,Y_N)$
means that there exist $a_1,\dots,a_N\in\R$ such that
$\sum_n a_n=0$ and $Y'_n=Y_n+a_n$.
We denote
$$
\EEE(X,Y)=\{\QQ\in\PPP:\exists Y'\sim Y\text{ such that }
(X,Y',\QQ)\text{ is an Arrow-Debreu-type equilibrium}\}.
$$

\begin{Theorem}
\label{EC7}
Assume that $M<\infty$
{\rm(}by Proposition~\ref{EC2}, this is equivalent to:
$\bigcap_n\DDD_n\cap\RRR(A_n)\ne\emp${\rm)}.
Let $X_n\in A_n$, $Y_n\in L_s^1(\DDD)$, $\sum_n Y_n=0$.
The following conditions are equivalent:
\begin{mitemize}
\item[\rm(i)] $(X,Y)$ is a Pareto-type equilibrium;
\item[\rm(ii)] there exist $\QQ\in\PPP$ and $Y'\sim Y$
such that $(X,Y',\QQ)$ is an Arrow-Debreu-type
equilibrium.
\end{mitemize}
If these conditions are satisfied, then
$$
\EEE(X,Y)=\EX_{\bigcap_n\!\DDD_n\cap\RRR(A_n)}(W).
$$

If each $A_n$ is a linear space, then {\rm(i)},
{\rm(ii)} are equivalent to:
\begin{mitemize}
\item[\rm(iii)] $\bigcap_n\EX_{\DDD_n}(W_n+X_n+Y_n)
\cap\RRR(A_n)\ne\emp$.
\end{mitemize}
Moreover, in this case
$$
\EEE(X,Y)=\bigcap_{n=1}^N\EX_{\DDD_n}(W_n+X_n+Y_n)\cap\RRR(A_n).
$$
\end{Theorem}

\sks
\Proof
\textit{Step~1.}
Let us prove the implication (i)$\Rightarrow$(ii).
Take $\QQ\in\EX_{\bigcap_n\!\DDD_n\cap\RRR(A_n)}(W)$
(this set is nonempty due to~\cite[Prop.~2.9]{C061}).
Using Proposition~\ref{EC2}, we can write
$$
\sum_{n=1}^N u_n(W_n+X_n+Y_n)
\le\sum_{n=1}^N\EE_\QQ(W_n+Y_n)
\le\EE_\QQ\sum_{n=1}^N(W_n+Y_n)
=\EE_\QQ W
=M
$$
(note that the expectation operator understood in the
sense of~\cite[Def.~2.3]{C061} has the property
$\EE(\xi+\eta)\ge\EE\xi+\EE\eta$).
Since the left-hand side and the right-hand side of the
above inequality coincide, we get
$u_n(W_n+X_n+Y_n)=\EE_\QQ(W_n+Y_n)$ for any~$n$.
We can find $Y'\sim Y$ such that $\EE_\QQ Y'_n=0$
for any~$n$.
Then, for any~$n$, $\xi\in A_n$, and $\eta\in L^1(\QQ)$
such that $\EE_\QQ\eta=0$, we have
$$
u_n(W_n+\xi+\eta)
\le\EE_\QQ(W_n+\xi+\eta)
\le\EE_\QQ W_n
=u_n(W_n+X_n+Y'_n).
$$
Thus, $(X,Y',\QQ)$ is an Arrow-Debreu-type equilibrium.

\sks
\textit{Step~2.}
Let us prove the implication (ii)$\Rightarrow$(i).
Suppose that there exist $\wt X_n\in A_n$ and
$\wt Y_n\in L_s^1(\DDD_n)$ with $\sum_n\wt Y_n=0$ such that
\begin{equation}
\label{ec4}
\sum_{n=1}^N u_n(W_n+\wt X_n+\wt Y_n)
>\sum_{n=1}^N u_n(W_n+X_n+Y_n).
\end{equation}

Fix $n$ and suppose that $\QQ\notin\DDD_n\cap\RRR(A_n)$.
By the Hahn-Banach theorem, there exists
$\eta\in L^\infty$ such that
$\EE_\QQ\eta<\inf_{\QQ\in\DDD_n\cap\RRR(A_n)}\EE_\QQ\eta$.
According to Proposition~\ref{GOG1}, there exists $\xi\in A_n$
such that $\EE_\QQ\eta<u_n(\xi+\eta)$.
This means that $u_n(\xi+\eta-\EE_\QQ\eta)>0$. Then
\begin{equation}
\label{ec5}
u_n(W_n+\al\xi+\al(\eta-\EE_\QQ\eta))\xra[\al\to\infty]{}\infty.
\end{equation}
In view of Proposition~\ref{EC2}, the condition $M<\infty$
implies that $\DDD_n\cap\RRR(A_n)\ne\emp$.
Then, for any $\wt\QQ\in\DDD_n\cap\RRR(A_n)$, we have
$$
u_n(W_n+X_n+Y_n)
\le\EE_{\wt\QQ}(W_n+Y_n)
<\infty,
$$
which contradicts~\eqref{ec5}.

Thus, $\QQ\in\bigcap_n\DDD_n\cap\RRR(A_n)$.
In particular, $Y_n\in L^1(\QQ)$, so that
we can find $\wt Y'\sim\wt Y$ such that $\EE_\QQ\wt Y'_n=0$
for any~$n$. Then
$$
\sum_{n=1}^N u_n(W_n+\wt X_n+\wt Y_n)
=\sum_{n=1}^N u_n(W_n+\wt X_n+\wt Y'_n)
\le\sum_{n=1}^N u_n(W_n+X_n+Y_n),
$$
which contradicts~\eqref{ec4}.

\sks
\textit{Step~3.}
It was shown in Step~1 that
$$
\EX_{\bigcap_n\!\DDD_n\cap\RRR(A_n)}(W)\subseteq\EEE(X,Y).
$$
Let us prove the reverse inclusion.
Take $\QQ\in\EEE(X,Y)$ and find $Y'\sim Y$ such that
$(X,Y',\QQ)$ is an Arrow-Debreu-type equilibrium.
It was shown in Step~2 that
$\QQ\in\bigcap_n\DDD_n\cap\RRR(A_n)$.
Applying Proposition~\ref{GOG1} to the $\DDD_n$-consistent
convex cone $A=\{\eta\in L^\infty:\EE_\QQ\eta=0\}$, we get
\begin{equation}
\label{ec6}
\sup_{\eta\in L^1(\QQ):\EE_\QQ\eta=0}u_n(W_n+\eta)
=\EE_\QQ W_n,
\quad n=1,\dots,N.
\end{equation}
Thus,
$$
M\ge\sum_{n=1}^N u_n(W_n+X_n+Y_n)
=\sum_{n=1}^N u_n(W_n+X_n+Y'_n)
\ge\sum_{n=1}^N\EE_\QQ W_n
=\EE_\QQ W.
$$
An application of Theorem~\ref{EC3} completes the proof.

\sks
\textit{Step~4.}
Assume that each $A_n$ is linear and let us prove
the inclusion
$$
\bigcap_{n=1}^N\EX_{\DDD_n}(W_n+X_n+Y_n)\cap\RRR(A_n)
\subseteq\EEE(X,Y).
$$
Take $\QQ$ from the left-hand side of this inclusion.
Find $Y'\sim Y$ such that $\EE_\QQ Y'_n=0$ for any~$n$.
Clearly, $\QQ\in\bigcap_n\EX_{\DDD_n}(W_n+X_n+Y'_n)\cap\RRR(A_n)$,
so that
$$
u_n(W_n+X_n+Y'_n)
=\EE_\QQ(W_n+X_n+Y'_n)
=\EE_\QQ(W_n+X_n),\quad
n=1,\dots,N.
$$
As $\EX_{\DDD_n}(W_n+X_n+Y'_n)\ne\emp$, we have
(by the definition of $\EX$) that
$u_n(W_n+X_n+Y'_n)>-\infty$.
Furthermore, (by the definition of $\RRR$)
$\EE_\QQ X_n\le0$, so that $X_n\in L^1(\QQ)$.
For any~$n$, $\xi\in A_n$, and $\eta\in L^1(\QQ)$ such that
$\EE_\QQ\eta=0$, we have
\begin{align*}
u_n(W_n+\xi+\eta)
&\le\EE_\QQ(W_n+\xi+\eta)
=\EE_\QQ(W_n+\xi)\\
&=\EE_\QQ(W_n+X_n+(\xi-X_n))
\le\EE_\QQ(W_n+X_n)
=u_n(W_n+X_n+Y'_n)
\end{align*}
(in the second inequality, we used the linearity of $A_n$),
so that $(X,Y',\QQ)$ is an Arrow-Debreu-type equilibrium.

\sks
\textit{Step~5.}
Assume that each $A_n$ is linear and let us prove the
inclusion
$$
\EEE(X,Y)\subseteq
\bigcap_{n=1}^N\EX_{\DDD_n}(W_n+X_n+Y_n)\cap\RRR(A_n).
$$
Take $\QQ\in\EEE(X,Y)$ and find $Y'\sim Y$ such that
$(X,Y',\QQ)$ is an Arrow-Debreu-type equilibrium.
It was shown in Step~2 that
$\QQ\in\bigcap_n\DDD_n\cap\RRR(A_n)$.
Applying~\eqref{ec6} and Proposition~\ref{GOG1}, we get
\begin{align*}
\EE_\QQ W_n
&=\sup_{\eta\in L^1(\QQ)\,:\,\EE_\QQ\eta=0}u_n(W_n+\eta)
\le u_n(W_n+X_n+Y'_n)\\
&\le\EE_\QQ(W_n+X_n+Y'_n)
\le\EE_\QQ W_n,\quad
n=1,\dots,N.
\end{align*}
Consequently, $u_n(W_n+X_n+Y'_n)=\EE_\QQ(W_n+X_n+Y'_n)$
for any~$n$, which means that
$\QQ\in\bigcap_n\EX_{\DDD_n}(W_n+X_n+Y'_n)
=\bigcap_n\EX_{\DDD_n}(W_n+X_n+Y_n)$.\End

\skm
The assumption that each $A_n$ is linear is essential
for the second part of Theorem~\ref{EC7} as shown by the
following example.

\begin{Example}\rm
\label{EC8}
Let $N=2$, $\DDD_1=\DDD_2=\{\PP\}$,
$A_1=A_2=\R_-$ (i.e. $A_1$, $A_2$ consist of
random variables that are identically equal to a
negative constant), and $W_1=W_2=0$.
Take $X_1=X_2=-1$, $Y_1=Y_2=0$.
Then $\bigcap_n\EX_{\DDD_n}(X_n+Y_n)\cap\RRR(A_n)=\{\PP\}$,
but clearly $(X,Y)$ is not a Pareto-type equilibrium.\End
\end{Example}

By Theorem~\ref{EC7}, the set $\EEE(X,Y)$ does not
depend on $(X,Y)$.
We call it the set of \textit{equilibrium measures}
and denote by~$\EEE$
(Theorem~\ref{EC7} yields the representation
$\EEE=\EX_{\bigcap_n\!\DDD_n\cap\RRR(A_n)}(W)$).

From the financial point of view, $\EEE$ is the set
of equilibrium price systems.
Thus, it is natural to define the set of
\textit{unconstrained equilibrium prices} of a contingent claim~$F$
simply as $I_E(F):=\{\EE_\QQ F:\QQ\in\EEE\}$.

%--------------------------------------------------------
\subsection{Constrained Equilibrium}
\label{EI}

Let $u_n$, $\DDD_n$, $A_n$, and $W_n$ be the same as in
the previous subsection.
Let $S$ be a $d$-dimensional random vector whose
components belong to $\bigcap_n L_s^1(\DDD_n)$.
From the financial point of view, there are $d$ financial
contracts that can be exchanged between the agents,
and $S^i$ means the payoff of the $i$-th contract.
(There is no relation between $S$ and $A_n$;
$A_n$ means the set of P\&Ls that can be obtained by
the $n$-th agent without trading the assets $1,\dots,d$.)

\begin{Definition}\rm
\label{EI1}
The \textit{maximal overall utility} is defined as
$$
M=\sup_{\vbox{
\vspace{1mm}
\hbox{\hspace*{6mm}$\scriptstyle X_n\in A_n$}
\vspace{-1mm}
\hbox{$\scriptstyle h_n\in\R^d\,:\,\sum_n\!h_n=0$}}}
\hspace{1mm}\sum_{n=1}^N u_n(W_n+X_n+\lb h_n,S\rb),
$$
where the sum is understood as $-\infty$ if any of the
summands equals $-\infty$.
\end{Definition}

Let us introduce the notation
\begin{align*}
G_n&=\cl\{\EE_\QQ S:\QQ\in\DDD_n\cap\RRR(A_n)\},\\
\wt G_n&=\cl\{\EE_\QQ(S,W_n):\QQ\in\DDD_n\cap\RRR(A_n)\},\\
f_n(x)&=\inf\{y:(x,y)\in\wt G_n\},\quad x\in G_n,\\
G&=\bigcap_{n=1}^N G_n,\\
f(x)&=\sum_{n=1}^N f_n(x),\quad x\in G.
\end{align*}

\begin{Proposition}
\label{EI2}
We have
$$
M=\inf_{x\in G}f(x),
$$
where $\inf\emp:=\infty$.
\end{Proposition}

\sksmin
\Proof
By Proposition~\ref{GOG1}, we have for any $h_n\in\R^d$,
\begin{align*}
\sup_{X_n\in A_n}u_n(W_n+X_n+\lb h_n,S\rb)
&=\inf_{\QQ\in\DDD_n\cap\RRR(A_n)}\EE_\QQ(W_n+\lb h_n,S\rb)\\
&=\inf_{x\in\wt G_n}\lb(h_n,1),x\rb\\
&=\inf_{x\in G_n}(\lb h_n,x\rb+f_n(x)).
\end{align*}
Standard results of convex analysis (see~\cite[Th.~16.4]{R97})
yield
$$
M=\sup_{h_n\in\R^d\,:\,\sum_n\!h_n=0}
\sum_{n=1}^N\inf_{x\in G_n}(\lb h_n,x\rb+f_n(x))
=\inf_{x\in G}f(x).
$$

We now pass on to various definitions of equilibrium.

\begin{Definition}\rm
\label{EI3}
A \textit{Pareto-type equilibrium} is a collection
$(X,h)=$ $(X_1,\dots,X_N,h_1,\dots,h_N)$ such that
\begin{mitemize}
\item[(a)] $X_n\in A_n$;
\item[(b)] $h_n\in\R^d$, $\sum_n h_n=0$;
\item[(c)] there do not exist $(X',h')$ satisfying (a), (b)
and such that
\begin{align*}
\forall n,\;u_n(W_n+X'_n+\lb h'_n,S\rb)\ge
u_n(W_n+X_n+\lb h_n,S\rb),\\
\exists n:u_n(W_n+X'_n+\lb h'_n,S\rb)
>u_n(W_n+X_n+\lb h_n,S\rb).
\end{align*}
\end{mitemize}
\end{Definition}

It is easy to see from the translation invariance
property ($u_n(X+m)=u_n(X)+m$) that condition~(c) is
equivalent to:
\begin{mitemize}
\item[(c')] $\sum_n u_n(W_n+X_n+\lb h_n,S\rb)=M$.
\end{mitemize}

\begin{Definition}\rm
\label{EI4}
An \textit{Arrow-Debreu-type equilibrium}
is a collection $(X,h,P)$, where $P\in\R^d$, such that
\begin{mitemize}
\item[(a)] $X_n\in A_n$;
\item[(b)] $h_n\in\R^d$, $\sum_n h_n=0$;
\item[(c)] for any $n$,
$$
u_n(W_n+X_n+\lb h_n,S-P\rb)
=\max_{\xi\in A_n,\,\eta\in\R^d}u_n(W_n+\xi+\lb\eta,S-P\rb).
$$
\end{mitemize}
\end{Definition}

Let us introduce the notation
$$
E(X,h)=\{P\in\R^d:(X,h,P)\text{ is an Arrow-Debreu-type
equilibrium}\}.
$$

\begin{Theorem}
\label{EI5}
Assume that $M<\infty$
{\rm(}by Proposition~\ref{EI2}, this is equivalent to:
$G\ne\emp${\rm)}.
Let $X_n\in A_n$, $h_n\in\R^d$, $\sum_n h_n=0$.
The following conditions are equivalent:
\begin{mitemize}
\item[\rm(i)] $(X,h)$ is a Pareto-type equilibrium;
\item[\rm(ii)] there exists $P\in\R^d$
such that $(X,h,P)$ is an Arrow-Debreu-type
equilibrium.
\end{mitemize}
If these conditions are satisfied, then
$$
E(X,h)=\argmin_{x\in G}f(x).
$$

If each $A_n$ is a linear space, then {\rm(i)},
{\rm(ii)} are equivalent to:
\begin{mitemize}
\item[\rm(iii)] $\bigcap_n\{\EE_\QQ S:
\QQ\in\EX_{\DDD_n}(W_n+X_n+\lb h_n,S\rb)\cap\RRR(A_n)\}\ne\emp$.
\end{mitemize}
Moreover, in this case
$$
E(X,h)=\bigcap_{n=1}^N\{\EE_\QQ S:
\QQ\in\EX_{\DDD_n}(W_n+X_n+\lb h_n,S\rb)\cap\RRR(A_n)\}.
$$
\end{Theorem}

\sksmin
\Proof
\textit{Step~1.}
Let us prove the implication (i)$\Rightarrow$(ii).
Take $P\in\argmin_{x\in G}f(x)$.
Using Proposition~\ref{EI2}, we can write
\begin{align*}
\sum_{n=1}^N u_n(W_n+X_n+\lb h_n,S\rb)
&\le\sum_{n=1}^N\inf_{x\in\wt G_n}\lb(h_n,1),x\rb\\
&\le\sum_{n=1}^N(f_n(P)+\lb h_n,P\rb)
=f(P)
=M.
\end{align*}
As the left-hand side and the right-hand side of this
inequality coincide, we get
$u_n(W_n+X_n+\lb h_n,S\rb)=f_n(P)+\lb h_n,P\rb$
for any~$n$. Thus, for any~$n$, $\xi\in A_n$, and
$\eta\in\R^d$, we have
\begin{align*}
u_n(W_n+\xi+\lb\eta,S-P\rb)
&=u_n(W_n+\xi+\lb\eta,S\rb)-\lb\eta,P\rb\\
&\le f_n(P)+\lb\eta,P\rb-\lb\eta,P\rb\\
&=u_n(W_n+X_n+\lb h_n,S-P\rb).
\end{align*}
Thus, $(X,h,P)$ is an Arrow-Debreu-type equilibrium.

\sks
\textit{Step~2.}
The implication (ii)$\Rightarrow$(i) follows from the
inequality: for any $\wt X_n\in A_n$ and
$\wt h_n\in\R^d$ with $\sum_n\wt h_n=0$, we have
\begin{align*}
\sum_{n=1}^N u_n(W_n+\wt X_n+\lb\wt h_n,S\rb)
&=\sum_{n=1}^N u_n(W_n+\wt X_n+\lb\wt h_n,S-P\rb)\\
&\le\sum_{n=1}^N u_n(W_n+X_n+\lb h_n,S-P\rb)\\
&=\sum_{n=1}^N u_n(W_n+X_n+\lb h_n,S\rb).
\end{align*}

\sks
\textit{Step~3.}
It was shown in Step~1 that
$$
\argmin_{x\in G}f(x)\subseteq E(X,h).
$$
Let us prove the reverse inclusion.
Take $P\in E(X,h)$.
Fix $n$ and suppose that $P\notin G_n$.
Then there exists $\eta\in\R^d$ such that
$\lb\eta,P\rb<\inf_{x\in G_n}\lb\eta,x\rb$.
According to Proposition~\ref{GOG1}, there exists $\xi\in A_n$
such that $\lb\eta,P\rb<u_n(\xi+\lb\eta,S\rb)$. Then
\begin{equation}
\label{ei1}
u_n(W_n+\al\xi+\lb\al\eta,S-P\rb)\xra[\al\to\infty]{}\infty.
\end{equation}
In view of Proposition~\ref{EI2}, the condition $M<\infty$
implies that $\DDD_n\cap\RRR(A_n)\ne\emp$.
Then, for any $\QQ\in\DDD_n\cap\RRR(A_n)$, we have
$$
u_n(W_n+X_n+\lb h_n,S-P\rb)
\le\EE_\QQ(W_n+\lb h_n,S-P\rb),
$$
which contradicts~\eqref{ei1}.

Thus, $P\in G$. Using Proposition~\ref{GOG1}, we can write
\begin{align*}
M&\ge\sum_{n=1}^N u_n(W_n+X_n+\lb h_n,S-P\rb)\\
&=\sum_{n=1}^N\sup_{\xi\in A_n,\eta\in\R^d}
u_n(W_n+\xi+\lb\eta,S-P\rb)\\
&=\sum_{n=1}^N\sup_{\eta\in\R^d}\inf_{x\in G_n}(f_n(x)+\lb\eta,x-P\rb)\\
&=\sum_{n=1}^N f_n(P)=f(P).
\end{align*}
An application of Proposition~\ref{EI2} completes the proof.

\sks
\textit{Step~4.}
Assume that each $A_n$ is linear and let us prove the
inclusion
$$
\bigcap_{n=1}^N\{\EE_\QQ S:
\QQ\in\EX_{\DDD_n}(W_n+X_n+\lb h_n,S\rb)\cap\RRR(A_n)\}
\subseteq E(X,h).
$$
Take $P$ from the left-hand side of this inclusion.
Using the same arguments as in the proof
of~\cite[Prop.~2.9]{C061}, we can find for every~$n$
a measure $\QQ_n\in\EX_{\DDD_n}(W_n+X_n+\lb h_n,S\rb)
\cap\RRR(A_n)$ such that $P=\EE_{\QQ_n}S$
(using the $\DDD_n$-consistency of $A_n$,
it is easy to check that
$\EX_{\DDD_n}(W_n+X_n+\lb h_n,S\rb)\cap\RRR(A_n)$
is $L^1$-closed, so this set is weakly compact). Then
\begin{align*}
u_n(W_n+X_n+\lb h_n,S-P\rb)
&=\EE_{\QQ_n}(W_n+X_n+\lb h_n,S\rb)-\lb h,P\rb\\
&=\EE_{\QQ_n}(W_n+X_n),\quad n=1,\dots,N.
\end{align*}
As $\EX_{\DDD_n}(W_n+X_n+\lb h_n,S\rb)\ne\emp$,
we have (by the definition of $\EX$) that
$u_n(W_n+X_n+\lb h_n,S\rb)>-\infty$.
Furthermore, (by the definition of $\RRR$) $\EE_{\QQ_n}X_n\le0$,
so that $X_n\in L^1(\QQ_n)$.
For any~$n$, $\xi\in A_n$, and $\eta\in\R^d$, we have
\begin{align*}
u_n(W_n+\xi+\lb\eta,S-P\rb)
&\le\EE_{\QQ_n}(W_n+\xi+\lb\eta,S-P\rb)
=\EE_{\QQ_n}(W_n+\xi)\\
&=\EE_{\QQ_n}(W_n+X_n+(\xi-X_n))
\le\EE_{\QQ_n}(W_n+X_n)\\
&=u_n(W_n+X_n+\lb h_n,S-P\rb)
\end{align*}
(in the second inequality, we used the linearity of $A_n$),
so that $(X,h,P)$ is an Arrow-Debreu-type equilibrium.

\sks
\textit{Step~5.}
Assume that each $A_n$ is linear and let us prove the
inclusion
$$
E(X,h)\subseteq\bigcap_{n=1}^N\{\EE_\QQ S:
\QQ\in\EX_{\DDD_n}(W_n+X_n+\lb h_n,S\rb)\cap\RRR(A_n)\}.
$$
Take $P\in E(X,h)$.
It was shown in Step~3 that $P\in G$.
Using the same arguments as in the proof
of~\cite[Prop.~2.9]{C061}, we can find for every~$n$
a measure $\QQ_n\in\DDD_n\cap\RRR(A_n)$
such that $\EE_{\QQ_n}S=P$ and $\EE_{\QQ_n}W_n=f_n(P)$.
Applying Proposition~\ref{GOG1}, we get
\begin{align*}
f_n(P)
&=\sup_{\eta\in\R^d}\inf_{x\in G_n}(f_n(x)+\lb\eta,x-P\rb)\\
&=\sup_{\xi\in A_n,\,\eta\in\R^d}u_n(W_n+\xi+\lb\eta,S-P\rb)\\
&=u_n(W_n+X_n+\lb h_n,S-P\rb)\\
&\le\EE_{\QQ_n}(W_n+X_n+\lb h_n,S-P\rb)\\
&\le\EE_{\QQ_n}W_n
=f_n(P).
\end{align*}
Consequently,
$u_n(W_n+X_n+\lb h_n,S-P\rb)=\EE_{\QQ_n}(W_n+X_n+\lb h_n,S-P\rb)$
for any~$n$, which means that
$\QQ_n\in\EX_{\DDD_n}(W_n+X_n+\lb h_n,S\rb)$.\End

\skm
Let $G_n^\circ$ denote the relative interior of $G_n$,
and, for $P\in G_n$, we denote
$$
N_{\wt G_n}(P)=\bigl\{\eta\in\R^d:\lb(\eta,1),(P,f_n(P))\rb
=\inf_{x\in\wt G_n}\lb(\eta,1),x\rb\bigr\}.
$$
If $\wt G_n$ has a nonempty interior, then $N_{\wt G_n}(P)$ is
the set of vectors $\eta\in\R^d$ such that $(\eta,1)$
is an inner normal to $\wt G_n$ at the point $(P,f_n(P))$.

\begin{figure}[!h]
\begin{picture}(150,92.5)(-80,-40)
\put(-40.5,-4.2){\includegraphics{poem.15}}
\put(-45,0){\vector(1,0){95}}
\put(-45,0){\vector(0,1){50}}
\multiput(-15,0.5)(0,1.95){10}{\line(0,1){1}}
\multiput(-10,0.5)(0,1.95){22}{\line(0,1){1}}
\multiput(-3,0.5)(0,2){15}{\line(0,1){1}}
\multiput(0,0.5)(0,1.95){17}{\line(0,1){1}}
\multiput(3,0.5)(0,2){10}{\line(0,1){1}}
\multiput(10,0.5)(0,1.95){22}{\line(0,1){1}}
\multiput(15,0.5)(0,1.9){5}{\line(0,1){1}}
\multiput(-44.5,32.8)(1.95,0){23}{\line(1,0){1}}
\multiput(-3,30)(2,0){2}{\line(1,0){1}}
\multiput(0,20)(2,0){2}{\line(1,0){1}}
\put(-2.5,14.7){\scalebox{0.85}{$\biggl\{$}}
\put(-0.8,25){\scalebox{0.85}{$\biggr\}$}}
\put(-3.2,14.7){\scalebox{0.8}{$1$}}
\put(1,25){\scalebox{0.8}{$1$}}
\put(-15,0){\circle*{1}}
\put(0,0){\circle*{1}}
\put(15,0){\circle*{1}}
\put(48,-3){\small $x$}
\put(-48,48.5){\small $y$}
\put(-32,16){\small $f_1(x)$}
\put(22,6){\small $f_2(x)$}
\put(-9,41){\scalebox{0.75}{$f_1(x)+f_2(x)$}}
\put(-27.5,35.3){\small $\wt G_1$}
\put(22.5,24.3){\small $\wt G_2$}
\put(-44,-7){\scalebox{0.8}{$I_{\text{\sl NBC}(\DDD_1,A_1,W_1)}(F)$}}
\put(-5,-8){\scalebox{0.8}{$\IE(F)$}}
\put(19.5,-7){\scalebox{0.8}{$I_{\text{\sl NBC}(\DDD_2,A_2,W_2)}(F)$}}
\put(-7.5,1){\small $h_1$}
\put(3.5,1){\small $h_2$}
\put(-50,32){\small $M$}
\put(-52.5,-27.5){\parbox{110mm}{\small\textbf{Figure~11.}
Geometric solution of the constrained equilibrium problem.
Here $S=F$ is a one-dimensional contingent claim.
The figure shows the maximal overall utility~$M$,
the equilibrium holdings $h_n$, and the
equilibrium price $I_E(F)$.
It also shows the NBC prices
$I_{\text{\sl NBC}(\DDD_n,A_n,W_n)}(F)$
of different agents.}}
\end{picture}
\end{figure}

\begin{Theorem}
\label{EI6}
Assume that $\bigcap_n G_n^\circ\ne\emp$.
Take $P\in\argmin_{x\in G}f(x)$.
Then there exist $h_n\in N_{\wt G_n}(P)$ such that
$\sum_n h_n=0$.
Assume that, for each~$n$, there exists
$X_n\in\argmax_{\xi\in A_n}u_n(W_n+\xi+\lb h_n,S\rb)$.
Then $(X_1,\dots,X_N,h_1,\dots,h_N,P)$ is an
Arrow-Debreu-type equilibrium.
Conversely, any Arrow-Debreu-type equilibrium has such
a form.
\end{Theorem}

\sksminus
\Proof
Denote $f_n^*(\eta)=\inf_{x\in G_n}(f_n(x)+\lb\eta,x\rb)$,
$\eta\in\R^d$.
Standard results of convex analysis (see~\cite[Th.~16.4]{R97})
guarantee that there exist $h_1,\dots,h_N\in\R^d$ such
that $\sum_n h_n=0$ and $\sum_n f_n^*(h_n)=f(P)$.
It follows from the line
$$
f(P)
=\sum_{n=1}^N f_n^*(h_n)
\le\sum_{n=1}^N(f_n(P)+\lb h,P\rb)
=f(P)
$$
that $h_n\in N_{\wt G_n}(P)$ for any~$n$.
By Proposition~\ref{GOG1},
\begin{align*}
u_n(W_n+X_n+\lb h_n,S\rb)
&=\sup_{\xi\in A_n}u_n(W_n+\xi+\lb h_n,S\rb)\\
&=\inf_{x\in\wt G_n}\lb(h_n,1),x\rb\\
&=f_n^*(h_n)
=f_n(P)+\lb h_n,P\rb,\quad n=1,\dots,N.
\end{align*}
Consequently, for any~$n$, $\xi\in A_n$, and $\eta\in\R^d$,
\begin{align*}
u_n(W_n+\xi+\lb\eta,S-P\rb)
&\le\inf_{x\in\wt G_n}\lb(\eta,1),x\rb-\lb\eta,P\rb\\
&=f_n^*(\eta)-\lb\eta,P\rb
\le f_n(P)\\
&=u_n(W_n+X_n+\lb h_n,S-P\rb),
\end{align*}
so that $(X,h,P)$ is an Arrow-Debreu-type equilibrium.

Conversely, let $(X,h,P)$ be an Arrow-Debreu-type
equilibrium.
According to Proposition~\ref{EI2}, $P\in\argmin_{x\in G}f(x)$.
Using the same arguments as in Step~1 of the proof of
Theorem~\ref{EI5}, we deduce that $h_n\in N_{\wt G_n}(P)$
for any~$n$.
The inclusion
$X_n\in\argmax_{\xi\in A_n}u_n(W_n+\xi+\lb h_n,S\rb)$
is clear from the definition of the Arrow-Debreu-type
equilibrium.\End

\skm
By Theorem~\ref{EI5}, the set $E(X,Y)$ does not
depend on $(X,Y)$.
We call it the set of \textit{equilibrium prices}
and denote by~$E$
(Theorem~\ref{EI5} yields the representation
$E=\argmin_{x\in G}f(x)$).

From the financial point of view, $E$ is the set
of equilibrium price vectors for the multidimensional
contract~$S$.
If $S=F$ is a one-dimensional contingent claim,
we call $E$ the set of \textit{constrained equilibrium
prices} of~$F$ and denote it by $I_E(F)$.

%=======================================================
\section{Conclusion}
\label{C}

In~\cite{C061} and in the present paper, we proposed
seven different pricing techniques.
They differ by the inputs they require and by the ideas
behind them.
These techniques are compared by Figure~12 and by Table~1.

\begin{picture}(150,110)(2,-55)
\put(0,-40){\includegraphics{poem.16}}
\put(52.5,-40){\small Utility-based NGD interval}
\put(50.8,-30){\small RAROC-based NGD interval}
\put(53.5,-20){\small Multi-agent NBC interval}
\put(50,-10){\small Agent-independent NBC price}
\put(55,41){\small Single-agent NBC prices}
\put(47,20){\small Unconstrained}
\put(50,16){\small equilibrium}
\put(55,12){\small price}
\put(80,20){\small Constrained}
\put(81,16){\small equilibrium}
\put(86,12){\small price}
\put(20,-55){\small\textbf{Figure~12.} The form of fair
prices provided by various techniques}
\end{picture}

\clearpage
\begin{figure}[p]
\tabcolsep=2mm
\begin{tabular}{|p{34.5mm}|p{22mm}|p{90mm}|}
\hline
\hfil\rule{0mm}{5mm}
\textbf{Pricing}\phantom{aaaaa}
\mbox{\phantom{aaaa}\textbf{technique}}\hfil&
\hfil\raisebox{-2mm}{\textbf{Inputs}}\hfil&
\hfil\raisebox{-2mm}{\textbf{Form of the price interval}}\hfil\\[2mm]
\hhline{|=|=|=|}
\rule{0mm}{5mm}Utility-based \mbox{No Good Deals}&
$\DDD$, $A$&
$\{\EE_\QQ F:\QQ\in\DDD\cap\RRR\}$\\[5.9mm]
\hline
\rule{0mm}{5mm}RAROC-based \mbox{No Good Deals}&
\mbox{$\PD$, $\RD$,} $A$, $R$&
\raisebox{-2mm}{$\Bigl\{\EE_\QQ F:\QQ\in\bigl(\frac{1}{1+R}\,\PD
+\frac{R}{1+R}\,\RD\bigr)\cap\RRR\Bigr\}$}\\[5.9mm]
\hline
Agent-independent \mbox{No Better Choice}&
$\PD$, $\RD$, $A$&
\rule{0mm}{6.5mm}${\Bigl\{\EE_\QQ F:\QQ\in\bigl(\frac{1}{1+R_*}\,\PD
+\frac{R_*}{1+R_*}\,\RD\bigr)\cap\RRR\Bigr\}}$,
\rule{0mm}{5mm}${\text{where }R_*=\sup_{X\in A}\RAROC(X)}$\\[2mm]
\hhline{|~|~|-|}
&&\rule{0mm}{6.5mm}${\Bigl\{\EE_\QQ F:\QQ\!\in\!\bigl(\frac{1}{1+R_*}\EX_\PD(X_*)
\!+\!\frac{R_*}{1+R_*}\EX_\RD(X_*)\bigr)\cap\RRR\Bigr\}}$,
\rule{0mm}{5mm}${\text{where }X_*\in\argmax_{X\in A}\RAROC(X)}$\\[2mm]
\hline
\rule{0mm}{5mm}Single-agent \mbox{No Better Choice}&
$\DDD$, $A$, $W$&
$\{\EE_\QQ F:\QQ\in\EX_\DDD(W)\cap\RRR\}$\\
\hhline{|~|~|-|}
&&\rule{0mm}{5mm}${\{\EE_\QQ F:\QQ\in\EX_\DDD(W)\}}$
provided that this is a singleton and
${u(W)=\max_{X\in A}u(W+X)}$\\[2mm]
\hline
\rule{0mm}{5mm}Multi-agent \mbox{No Better Choice}&
${\DDD_1,\dots,\DDD_N}$, $A$, ${W_1,\dots,W_N}$&
$\bigl\{\EE_\QQ F:\QQ
\in\conv_{n=1}^N(\EX_{\DDD_n}(W_n)\cap\RRR)\bigr\}$\\[11mm]
\hline
\rule{0mm}{5mm}\mbox{Unconstrained} \mbox{equilibrium}&
$\DDD_1,\dots,\DDD_N$, $A_1,\dots,A_N$, $W_1,\dots,W_N$&
${\{\EE_\QQ F:\QQ\in\EEE\}\text{, where}}$
\rule{0mm}{5mm}${\EEE=\EX_{\bigcap_n\!\DDD_n\cap\RRR(A_n)}\bigl(\sum_n\!W_n\bigr)}$\\[5mm]
\hline
\rule{0mm}{5.5mm}\mbox{Constrained} \mbox{equilibrium}&
$\DDD_1,\dots,\DDD_N$, $A_1,\dots,A_N$, $W_1,\dots,W_N$&
${\argmin_x\sum_n f_n(x),\text{ where}}$
\rule{0mm}{5mm}${f_n(x)=\inf\{\EE_\QQ W_n:\QQ\in\DDD_n\cap\RRR(A_n),\,
\EE_\QQ F=x\}}$\\[4mm]
\hline
\end{tabular}

\vspace{7mm}
\hspace*{20mm}\parbox{120mm}{\small\textbf{Table~1.}
The form of fair prices provided by various techniques}
\end{figure}

%========================================================


\clearpage
\begin{thebibliography}{100}
\label{R}

\bibitem{A04} \textit{C.~Acerbi.}
Coherent representations of subjective risk aversion.
In: G.~Szeg\"o (Ed.). Risk measures for the 21st century.
Wiley, 2004, pp.~147--207.

\bibitem{BE05} \textit{P.~Barrieu, N.~El~Karoui.}
Inf-convolution of risk measures and optimal risk transfer.
Finance and Stochastics, {\bf 9} (2005), p.~269--298.

\bibitem{B95} \textit{F.~Black.}
Estimating expected return.
Financial Analysts Journal, 1995, p.~168--171.

\bibitem{CGM01} \textit{P.~Carr, H.~Geman, D.~Madan.}
Pricing and hedging in incomplete markets.
Journal of Financial Economics, {\bf 62} (2001),
p.~131--167.

\bibitem{C061} \textit{A.S.~Cherny.}
Pricing with coherent risk.
Preprint, available at:
\texttt{http://mech.math.msu.su/\~{}cherny}.

\bibitem{C05e} \textit{A.S.~Cherny.}
Weighted V@R and its properties.
To be published in Finance and Stochastics. Available at:
\texttt{http://mech.math.msu.su/\~{}cherny}.

\bibitem{CK99} \textit{J.~Cvitani\'c, I.~Karatzas.}
On dynamic measures of risk.
Finance and Stochastics, {\bf 3} (1999),
\mbox{p.~451--482}.

\bibitem{D05} \textit{F.~Delbaen.}
Coherent monetary utility functions.
Preprint, available at
\texttt{http://www.math.ethz.ch/\~{}delbaen}
under the name ``Pisa lecture notes''.

\bibitem{F03} \textit{T.~Fischer.}
Risk capital allocation by coherent risk measures based
on one-sided moments.
Insurance: Mathematics and Economics {\bf 32} (2003),
No.~1, p.~135--146.

\bibitem{HK04} \textit{D.~Heath, H.~Ku.}
Pareto equilibria with coherent measures of risk.
Mathematical Finance, {\bf 14} (2004), p.~163--172.

\bibitem{JPR05} \textit{A.~Jobert, A.~Platania,
L.C.G.~Rogers.}
A Bayesian solution to the equity premium puzzle.
Preprint, available at:
\texttt{http://www.statslab.cam.ac.uk/\~{}chris}.

\bibitem{JST05} \textit{E.~Jouini, W.~Schachermayer, N.~Touzi.}
Optimal risk sharing for law invariant monetary utility functions.
Preprint, available at:
\texttt{http://www.fam.tuwien.ac.at/\~{}wschach/pubs}.

\bibitem{K90} \textit{T.C.~Koopmans.}
Three essays on the state of economic science.
A.M.~Kelley, 1990.

\bibitem{L97} \textit{E.L.~Lehmann.}
Testing statistical hypotheses.
Springer, 1997.

\bibitem{M59} \textit{H.~Markowitz.}
Portfolio selection. Wiley, 1959.

\bibitem{N03} \textit{Y.~Nakano.}
Minimizing coherent risk measures of shortfall in
discrete-time models under cone constraints.
Applied Mathematical Finance, {\bf 10} (2003),
p.~163--181.

\bibitem{N04} \textit{Y.~Nakano.}
Minimization of shortfall risk in a jump-diffusion model.
Statistics and Probability Letters, {\bf 67} (2004),
p.~87--95.

\bibitem{R97} \textit{R.T.~Rockafellar.}
Convex analysis. Princeton, 1997.

\bibitem{RU00} \textit{R.T.~Rockafellar, S.~Uryasev.}
Optimization of conditional Value-At-Risk.
Journal of Risk, {\bf 2} (2000), No.~3, p.~21-41.

\bibitem{RUZ05} \textit{R.T.~Rockafellar, S.~Uryasev,
M.~Zabarankin.}
Master funds in portfolio analysis with general
deviation measures.
Journal of Banking and Finance, {\bf 29} (2005).

\bibitem{Sek04} \textit{J.~Sekine.}
Dynamic minimization of worst conditional expectation
of shortfall.
Mathematical Finance, {\bf 14} (2004), No.~4, p.~605--618.

\bibitem{S64} \textit{W.~Sharpe.}
Capital asset prices: A theory of market equilibrium
under conditions of risk.
Journal of Finance, {\bf 19} (1964), 425-442.

\end{thebibliography}
\end{document}